\documentclass[12pt]{amsart} 
\newfont{\sheaf}{eusm10 scaled\magstep1}

\newcommand{\C}{\ensuremath{\mathbb{C}}}
\newcommand{\R}{\ensuremath{\mathbb{R}}} 
\newcommand{\Z}{\ensuremath{\mathbb{Z}}} 
 
\newcommand{\hol}{\ensuremath{\mathcal{O}}} 
\newcommand{\PP}{\ensuremath{\mathbb{P}}} 
\newcommand{\Proof}{{\it Proof. }}

\newtheorem{teo}{Theorem}[section] 
\newtheorem{df}[teo]{Definition} 
\newtheorem{lem}[teo]{Lemma} 
\newtheorem{cor}[teo]{Corollary} 
\newtheorem{oss}[teo]{Remark} 
\newtheorem{prop}[teo]{Proposition}

\title{Deformation types of real and complex manifolds}
\author{Fabrizio Catanese}

\thanks{
The present research took place in the framework of the Schwerpunkt
"Globale Methode in der komplexen Geometrie", and of EAGER.
}
\begin{document} 
\thispagestyle{empty}

\maketitle

The present article grew  directly out of the lecture
delivered at the Chen-Chow Symposium, entitled "Algebraic surfaces : real
structures, topological and differentiable types", and  partly evolved
through the transparencies I used   later on to lecture at FSU
Tallahassee, and at  Conferences in Marburg and Napoli.
 
 The main reason however to change the title was that later on,
motivated by some open problems I mentioned in the Conference, I started
 to include in the article new results in higher dimensions, and then
  the theme of deformations in the large of real and complex manifolds
emerged as a central one.

For instance, we could summarize the Leit-Faden of the article as a
negative answer to the question whether, for a compact complex manifold
which is a $K( \pi , 1)$ (we consider  in some way these manifolds as
being simple objects),  the diffeomorphism type determines the
deformation type.

The first counterexamples to this question go essentially back to some old
papers of Blanchard and of Calabi (cf.  \cite{bla1},\cite{bla2}, \cite{cal}
\cite{somm}), showing the existence of non K\"ahler complex structures on
the product of a curve  with a two dimensional complex torus (the curve has
odd genus $g \geq 3$ for Calabi's examples which are, I believe, a special
case of the ones obtained by Sommese via the Blanchard method). Our
contribution here, in the more technical section 4, is to show that any
deformation in the large of a product of a curve of genus $g \geq 2$ with a
complex torus is again a variety of the same type. We show also that the
Sommese-Blanchard complex structures  give rise to infinitely
many Kuranishi families whose dimension is unbounded , but we are yet
unable to decide whether there are indeed infinitely many distinct
deformation types.

As a preliminary result , we  give a positive answer to a question
raised by Kodaira and Spencer (cf. \cite {k-s}, Problem 8, section 22,
page 907 of volume II of Kodaira's collected works), showing that
any deformation in the large of a complex torus is again a complex torus.
In the same section we also  give a criterion for a complex manifold to be a
complex torus, namely, to have the same integral cohomology algebra of a
complex torus, and to possess $n$ independent holomorphic and d-closed
1-forms. 

Whether complex tori admit other complex structures  with trivial
canonical bundle is still an open question (the Sommese-Blanchard
examples  always give, on 3- manifolds diffeomorphic to tori, complex
structures with effective anticanonical bundle).

One might believe that the previous "pathologies"  occur because we did not
 restrict ourselves to the class of  K\"ahler manifolds. However, even the
K\"ahler condition is not sufficient to ensure that deformation and
diffeomorphism type coincide.
 
 In fact, and it is striking that we have (non rigid) examples already in
dimension $2$, we illustrate in the last sections
(cf.\cite{cat7},\cite{cat8} for full proofs), that there are  K\"ahler
surfaces which are $K( \pi , 1)$'s, and for which there are different
deformation types, for a fixed differentiable structure.

These final examples tie up perfectly with the beginning of our journey,
that is, the   classification problem of real varieties.
 
  The original point of view illustrated in the lecture was to 
 zoom  the focus , considering  wider and wider classes of
objects:

{\bf  CLASSES OF OBJECTS }

\begin{itemize}

\item Smooth real algebraic Varieties/ Deformation
\item Smooth complex algebraic Varieties/ Deformation
\item Complex K\"ahler Manifolds/ Deformation
\item Symplectic Manifolds
\item Differentiable Manifolds
\item Topological Manifolds
\end{itemize}

In the article we could not illustrate all the possible aspects of the
problem ( we omitted for instance the topic of symplectic manifolds, see
\cite{au} for a nice survey), we tried however  to  cover most of them.

The paper is organized as follows: in section 1 we recall the definition
and the main properties and examples of real varieties, and we report on
recent results (\cite{c-f2}) on the Enriques classification of real
surfaces.

Section 2 is devoted to some open questions, and illustrates the simplest
classification problem for real varieties, namely the case of  real curves
of genus one.

In the next section we illustrate  the role of hyperelliptic surfaces in
the Enriques classification, and explain the main techniques used in the
classification of the real hyperelliptic surfaces (done in \cite{c-f2}).

Section 4  is mainly devoted, as already
mentioned, to complex tori and to  products  $C \times T$, and there we
barely mention the constructions of  Sommese, Blanchard and Calabi. We end
up with the classification and deformation theory of real tori (this is
also a new result).

Section 5 is devoted to a detailed description of what we call the
Blanchard-Calabi 3-folds ( they are obtained starting from curves in the
Grassmann variety $G(1,3)$, via a generalization of the previous
constructions). This approach allows to construct smooth
submanifolds of the Kuranishi family of deformations of a Blanchard-Calabi
3-fold, corresponding to the deformations of
$X$ which preserve the  fibration with fibres 2-dimensional tori. 

We are moreover able to show that these submanifolds coincide with the
Kuranishi family for the Sommese-Blanchard 3-folds, and more generally for
the Blanchard-Calabi 3-folds whose associated ruled surface is non
developable.

In this way we obtain infinitely many families whose dimension,
for a fixed differentiable structure, tends to infinity. Finally, we sketch
the Calabi construction of almost complex structures, without proving in
detail that the Calabi examples are a special case of the Sommese-Blanchard
3-folds.

In section 6 we recall the topological classification of simply connected
algebraic surfaces and briefly mention recent counterexamples to the
Friedman-Morgan speculation that for surfaces of general type the
deformation type should be determined by the differentiable type.

Sections 7 is devoted to some results and examples of triangle curves,
whereas the last section explains some details of the construction of
our counterexamples, obtained by suitable quotients of products of curves.
For these, the moduli space for a fixed differentiable type is the same
as the moduli space for a fixed topological type (and even the moduli space
for a given fundamental group and Euler number, actually!), and it has two
connected components, exchanged by complex conjugation.

\section{What is a  real variety?}

Let's then start with the first class, explaining
 what is a real variety, and what are the problems one is interested in
 (cf. \cite{dekha2} for a broader survey).

In general, e.g. in real life, one wants to solve
polynomial equations with real coefficients, and find real
solutions. Some theory is needed for this.

First of all, given a system of polynomial equations, in order to have some
continuity of the dependence of the solutions upon the choice
of the coefficients, one has to reduce it to a system of
homogeneous equations

$$ f_1(z_0, z_1, ... z_n) = 0 \\$$
$$........................$$
$$f_r(z_0, z_1, ... z_n) = 0.$$

The set $X$ of non trivial complex solutions of this system is
called an algebraic set of the projective space $\PP^n_{\C}$, and
the set $X(\R)$ of real solutions will be the intersection
$ X \cap \PP^n_\R$.

One notices that the set  $\bar{X}$ of complex conjugate
points is also an algebraic set, corresponding to the system
where we take the polynomials 
$$ \bar{f_j}(\bar{z_0},  ... \bar{z_n}) = 0 \\$$
( i.e., where we  conjugate the coefficients of the $f_j's$).

$X$ is said to be $\bf{real}$ if we may assume the $f_j's$
to have real coefficients, and in this case $X= \bar{X}$.

Assume that $X$ is real: for reasons stemming from Lefschstz' s
topological investigations, it is better to look at $X(\R)$ as
the subset of $X$ of the points which are left fixed by
the self mapping $\sigma$ given by complex conjugation
(curiously enough, Andre' Weil showed that one should use a
similar idea to study equations over finite fields, in this case
one lets $\sigma$ be the map raising each variable $z_i$ to
its $q$-th power ).

The simplest example of a real projective variety (a variety
is an algebraic set which is not the union of two proper
algebraic subsets) is a plane conic $ C \subset \PP^2_\C$.
It is defined by a single quadratic equation, and if it is
smooth, after a linear change of variables, its equation can
be written ( possibly multiplying it by $-1$) as

$$z_0^2 + z_1^2 + z_2^2 = 0 $$  or as 
$$z_0^2 - z_1^2 - z_2^2 = 0. $$

In the first case $X(\R) = \emptyset$, in the second case we
have that $X(\R)$ is a circle.

Indeed, the above cases exhaust the classification of smooth
real curves of genus $=0$.

We encounter genus $=1$ when we proceed to an equation of
degree $3$.

For instance, if we consider the family of real cubic
curves with affine equation 

$ y^2 - x^2 (x +1) = t$, we obtain 

\begin{itemize}
\item
two ovals for $ t < 0$
\item
one oval for $ t > 0$.
\end{itemize}

However, every polynomial of  odd degree has a real root,
whence every real cubic curve has $X(\R) \neq \emptyset$.

On the other hand, there are curves of genus $1$ and without
real points, as we were taught by Felix Klein (\cite{kl}) and his famous 
$\bf Klein's \ bottle$.

Simply look at a curve of genus $1$ as a complex torus
$X= \C / \Z + \tau \Z$. For instance, take $\tau = i$
and let $\sigma$ be the antiholomorphic self map induced by
$\sigma (z) = \bar{z} + 1/2$. By looking at the real part, we
see that there are no fixpoints, thus $X(\R) = \emptyset$:
indeed the well known Klein bottle is  exactly the quotient $X / \sigma$,
as an easy picture shows.
[Remark: this real curve  $X$ is the locus of
zeros of two real quadratic polynomials in $ \PP^3_\C$.] 

The previous example shows once more the usefulness of the
notion of an abstract manifold ( or variety), which is
indeed one of the keypoints in the classification theory.

\begin{df} A smooth real variety is a pair $(X, \sigma)$ 
 consisting of  a smooth complex variety $X$
of  complex dimension $n$
and of an  antiholomorphic involution $\sigma: X \rightarrow
X$ (involution means: $\sigma^2= Id$).
\end{df}

I.e., let $M$ be the differentiable manifold  underlying  $X$ 
and let $J$ be the complex structure of $X$ ($J$ is the
linear map on real tangent vectors provided  by
multiplication by $i$): then  the complex structure
$-J$ determines a complex manifold  which is called the
conjugate of $X$ and denoted by
$\bar{X}$, and  $\sigma$
is said to be antiholomorphic if it provides a 
holomorphic map between
the complex manifolds $X$ and $\bar{X}$.

{\bf Main problems}
(let's assume $X$ is compact)
\begin{itemize}
\item
Describe  the isomorphism classes
of such pairs $(X, \sigma)$.
\item
Or, at least describe the possible topological or differentiable types
of the pair $(X, \sigma)$.
\item
At least describe the possible topological types for the real part
$X': = X({\R}) = Fix(\sigma)$.
\end{itemize}

\begin{oss}

One can generalize the last problem and consider real pairs
$( Z \subset X , \sigma)$. 

Recall  indeed that Hilbert's 16-th problem is a special case
of  the last question  for the special case where $X= \PP^2$,
and $Z$ is a  smooth curve. In practice, the problem consists 
then in finding how many ovals $Z(\R)$ can have, and what is
their mutual disposition ( one inside the other , or not).
\end {oss}

For real algebraic curves one has the following nice 

 HARNACK's INEQUALITY

\begin{teo} Let $(C, \sigma)$ be a real curve of genus $g$.
Then the real part $C(\R)$ consists of a disjoint union of
$t$ ovals = circles ( topologically: $S^1$), where $ 0 \leq t
\leq g+1$.
\end{teo}

\begin{oss}
Curves with $g+1$ ovals are called M(aximal)-CURVES. 
An easy example of M-curves is provided by the hyperelliptic curves:
take  a real polynomial $P_{2g+2} (x)  $ of degree $2g + 2$
and with all the roots real, and consider the hyperelliptic
curves with affine equation 

$$ z^2= P_{2g+2} (x)  .$$ 
\end{oss}

The study of real algebraic curves has been the object of many
deep investigations, and although many questions remain still
open, one has a good knowledge of them, for instance Sepp\"ala
and Silhol (\cite{se-si}) and later Frediani (\cite{fre}) proved  

\begin{teo}
Given an integer $g$, the subset of the moduli space of complex curves
of genus $g$,  given by the curves which admit a real structure,
 is connected.
\end{teo}

A clue point to understand the meaning of the above theorem is that a
complex variety can have  several real structures, or none.
In fact the group $Aut(X)$ of biholomorphic automorphisms sits as a
subgroup of index at most $2$ in the group $Dian(X)$ consisting of 
$Aut(X)$ and of the antiholomorphic automorphisms. The real structures
are precisely the elements of order $2$ in $Dian(X) - Aut(X)$.
Therefore, as we shall see, the map of a  moduli space of real varieties
to the real part of the moduli space of complex structures can have
positive degree over some points, in particular, although the moduli space
of real curves of genus $g$ is not connected, yet its image in the moduli
space of complex curves is connected!

In higher dimensions, there are many fascinating questions,
and the next natural step is the investigation of the case of
algebraic surfaces, which is deeply linked to the
intriguing mystery of smooth 4- manifolds.

Indeed, for complex projective surfaces, we have the Enriques'
classification of surfaces up to birational equivalence
( equivalently, we have the classification of minimal surfaces,
i.e., of those surfaces
$S$ such that any holomorphic map $S \rightarrow S'$ of degree
$1$ is an isomorphism).

The Enriques (-Kodaira) classification of algebraic
varieties should consist in subdividing  the varieties $X$
according to their so-called Kodaira dimension, which is   a
number $kod(X) \in \{ - \infty, 0, 1, .. , dim(X)\}$ , and then
giving a detailed description of the varieties of $\bf special
\ type$, those for which  $kod(X) < dim (X)$ ( the varieties for
which $kod(X) = dim (X)$ are called of $\bf general \ 
type$).

As such, it has  been achieved for curves:

$$
\begin{tabular}{|c|c||c|}
\hline Curve & $Kod$  & $g$  \\
\hline
$\PP^1_{\C}$ & $- \infty$ & $0$   \\
Elliptic curve: $\C / \Z + \tau \Z $ & $0$ & $1$  \\
Curve of general type & $1$ & $ \geq 2$\\
\hline
\end{tabular}
$$ 
and for surfaces, for instance the following is 
 the Enriques'
classification of complex projective surfaces (where $C_g$ stands for a
curve of genus $g$):

$$
\begin{tabular}{|c|c||c||c|}
\hline Surface & $Kod$  & $P_{12}$ & Structure \\
\hline
 Ruled surface with  $q=g$& $- \infty$ & $0$&
$\PP^1_{\C} \times C_g$  \\ 
\hline
Complex torus 
& $0$ & $1$ & $\C^2 / \Lambda_4 $\\
\hline
K3 Surface & $0$& $1$&homeo to $X^2_4 \subset \PP^3_{\C}$ \\

\hline
 Enriques surface & $0$& $1$& $ K3/ (\Z /2 )$ \\
\hline
Hyperelliptic surface 
& $0$& $1$&  $(C_1 \times C'_1) / G$\\
\hline
Properly elliptic surface& $1$& $\geq 2$& $ dim \ \phi_{12}(S)=1$
, ...\\
\hline
Surface of general type & $2$
& $ \geq 2$ & $ dim \ \phi_{12}(S)=2$, ?\\

\end{tabular}
$$

Now, the Enriques classification of real algebraic
surfaces has not yet been achieved  in its strongest
form,  however  it has been achieved for Kodaira dimension
$0$ , thanks to  Comessatti (around $1911$, cf.
\cite{com1},\cite{com2},\cite{com3}), Silhol, Nikulin, Kharlamov and
Degtyarev (cf. references) and was finished in our joint work with Paola
Frediani (\cite{c-f2}) where we proved:

\begin{teo} Let $(S, \sigma)$ be a real hyperelliptic surface. 

1) Then the
differentiable type of the pair $(S, \sigma)$ is completely determined 
by the orbifold fundamental group exact sequence.\\
2) Fix the topological type of $(S, \sigma)$ corresponding to a
real hyperelliptic surface. Then the moduli space of the real surfaces $(S',
\sigma ')$ with the given topological type is
irreducible (and connected).\\
3) Real hyperelliptic surfaces fall into  $78$
topological types. In particular, the real part 
 $S(\R))$ of a real hyperelliptic surface is either
 \begin{itemize}
 \item
 a disjoint union of $c$ tori, where $ 0 \leq c \leq 4$
 \item
 a disjoint union of $b$ Klein bottles,
 where  $ 1 \leq b \leq 4$.
 \item
 the disjoint union of a torus and of a Klein bottle
 \item
 the disjoint union of a torus and of two Klein bottles.
 \end{itemize} 
\end{teo}
This result confirms a conjecture of Kharlamov that 
for real surfaces of Kodaira dimension $\leq 1$ the deformation
type of $(S, \sigma)$ is determined by the topological
type of $(S, \sigma)$.
Our method consists in 

1) Rerunning the classification theorem with special attention to
the real involution.

2) Finding out the primary role of the 
{\bf orbifold fundamental group},
which is so defined:

For $X$ real smooth, we have a double covering
$\pi : X \rightarrow Y = X/<\sigma>$ ($Y$ is called the
$\bf Klein \ variety$ of
$(X, \sigma)$),  ramified on the   real part of
$X$,
$ X({\R}) = Fix(\sigma)$.

In the case where  $X({\R}) = \emptyset $, 
 $ \pi_1^{orb}(Y)$ is just defined as
 the fundamental group $ \pi_1(Y)$.
 
 Otherwise, pick a point $x_0 \in X({\R})$, thus 
 $ \sigma (x_0) = x_0$ and the action of $\sigma$  on 
 $\pi_1(X, x_0)$ allows to 
  define $ \pi_1^{orb}(Y)$
 as a semidirect product of
 $\pi_1(X, x_0)$ with $\bf Z /2$.

One checks that the definition is independent of the choice of 
$x_0$. 

We have thus in all cases an exact sequence

$ 1 \rightarrow \pi_1(X) \rightarrow \pi_1^{orb}(Y) \rightarrow {\Z}/2
\rightarrow 1 $.

This sequence is very important when  $X$ is a $ K( \pi,1)$, i.e., 
when the
universal cover of $X$ is contractible.

\section{Complex manifolds which are $ K( \pi,1)$'s and real 
curves of genus 1. }

Some interesting questions are, both for complex and real
varieties

{\bf Question 1 : to what extent,  if $X$ is a $ K( \pi,1)$, does then
$\pi_1(X)$ also determine the differentiable type of $X$, not
only  its homotopy type? }

{\bf Question 2 : analogously, if $(X, \sigma)$ is a $ K( \pi,1)$ and is
real,  how much from the differentiable viewpoint is $(X, \sigma)$
determined by the  orbifold fundamental group sequence?}

There are for instance, beyond the case of hyperelliptic surfaces, other
similar instances in the Kodaira classification of real surfaces  
(joint work  in progress with Paola Frediani).

{\bf Question 3 : determine, for the real varieties  whose 
differentiable type
is determined by the orbifold  fundamental group, 
those for which the
 moduli spaces are irreducible and connected.}
 
 \begin{oss}
 However, already   for complex
manifolds,  the question whether (fixed the differentiable structure) there
is a unique deformation type,
 has several negative answers, as we
shall see in the later sections.
\end{oss}

 To explain   the seemingly  superfluous statement "irreducible
and connected", let us observe that : a hyperbola
$
\{ (x,y)
\in {\bf R^2} | xy =1
\}
$
 is irreducible but not connected. 
 
 I want to show now an easy example, leading to the quotient of
the above
 hyperbola by the involution $ (x,y) \rightarrow (-x,-y)$
 as a moduli
 space, and explaining the basic philosophy underlying the
 two  above  mentioned theorems concerning the orbifold fundamental group.
 
  REAL  CURVES OF GENUS ONE
 
 Classically,  the topological type of real curves  of genus $1$  is classified according to the number
$\nu = 0,1, or 2$ of  connected components
(= ovals = homeomorphic to circles) of their real part . By
abuse of language we shall also say: real elliptic curves, instead of
curves of genus $1$, although for many authors an elliptic curve comes
provided with one point defined over the base field (viz. : the origin!).

The orbifold fundamental group sequence is in this case 

$1 \rightarrow H_1(C, \Z) \cong \Z^2 \ \rightarrow 
\pi_1^{orb}(C) \rightarrow {\Z}/2 \rightarrow 1$.

and it splits iff $C(\bf R) \neq \emptyset$
(since $\pi_1^{orb}(C)$ has 
a representation as a group of affine 
transformations of the plane).

Step 1 : If there are no fixed points, the action $s$ of ${\Z}/2$ on 
$\Z^2$  (given by conjugation) is diagonalizable.

PROOF:

let $\sigma$ be represented by the affine transformation
$(x,y) \rightarrow s(x,y) + (a,b) $.
Now, $s$ is not diagonalizable
if and only if for a suitable  basis choice, $s(x,y)=
(y,x)$.  In this case the square of $\sigma$ is the
transformation $(x,y) \rightarrow (x,y) + (a+b,a+b) $, thus
$a+b$ is an integer, and therefore the points $(x, x-a)$
yield a fixed $S^1$ on the elliptic curve.
\hfill Q.E.D. \\

Step 2 : If there are no fixed points, moreover,  the
translation vector of the affine transformation inducing
$\sigma$ can be chosen to be
$1/2$ of the $+1$-eigenvector $e_1$ of $s$.
Thus we have exactly one normal form.

Step 3 :  there are exactly $3$ normal forms, and they are distinguished
by the values $0,1,2$ for  $\nu$. Moreover, 
 $s$ is diagonalizable if and only if $\nu$ is even.

PROOF: if there are fixed points, then $\sigma$ may be assumed
to be linear, so there are exactly two normal forms, according
to whether $\sigma$ is diagonalizable or not.
One sees immediately that  $\nu$ takes
then the respective values $2,1$.! Moreover, we have then verified that 
 $s$ is diagonalizable if and only if $\nu$ is even.
\hfill Q.E.D. \\

We content ourselves now with

THE DESCRIPTION OF THE MODULI SPACE FOR THE  
 CASE  $\nu =1$.

 $\sigma $ acts as follows: $ (x,y) \rightarrow (y,x) $.
 
We look then for  the complex structures $J$ which make $\sigma $
antiholomorphic,
i.e., we seek for the matrices $J$ with $J^2 = -1$ and with $ J s = - s J$.

The latter condition singles out the matrices
$\left(\begin{array}{cc}
a & b \\
-b & -a
\end{array} \right)$

while the first condition is equivalent to requiring that the
characteristic polynomial be equal to $ \lambda^2 + 1$,
whence, equivalently, $b^2 - a^2 =1$.

We get a hyperbola with two branches which are exchanged under
the involution $ J \rightarrow -J$, but, as we already remarked,
$J$ and $-J$ yield isomorphic real elliptic curves, thus
the moduli space is irreducible and connected.
\hfill Q.E.D. \\

\section{Real hyperelliptic surfaces on the scene}

The treatment of real hyperelliptic surfaces is similar
to the above sketched one of elliptic curves, and very much related to it,
because the orbifold group has an affine representation.
Recall:
\begin{df} A complex surface $S$ is said to be hyperelliptic 
if $S \cong (E
\times F)/G$, where $E$ and $F$ are elliptic curves and $G$
is a finite
group of
translations of $E$ with a faithful action on $F$ such that
$F/G \cong {\bf
P}^1$.
\end{df}

HISTORICAL REMARK : the points of elliptic curves
 can be parametrized by meromorphic functions of $z\in\C$.
 Around 1880, thanks to
the work of Appell, Humbert and Picard, there was much
interest for the "hyperelliptic varieties" of higher dimension
$n$, whose points
can be parametrized by meromorphic functions of $z \in \C^n$,
but not by rational functions of $z \in \C^n$. 
The classification of
hyperelliptic surfaces was finished by Bagnera and De Franchis
(\cite{b-df1},\cite{b-df2})
who were awarded the Bordin Prize in 1908 for this important
achievement.

\begin{teo}
\label{bdf} (Bagnera - de Franchis) Every hyperelliptic surface is one of the
following, where $E$, $F$ are  elliptic curves and $G$ is a group of
translations
of $E$ acting on $F$ as specified ($\rho$ is a primitive third root of
unity):
\begin{enumerate}
\item $(E \times F)/G$, $G = {\Z}/2$ acts on $F$ by $x
\mapsto -x$.
\item $(E \times F)/G$, $G = {\Z}/2 \oplus {\Z}/2$ acts on $F$ by $x
\mapsto -x$, $x \mapsto x +\epsilon$, where
$\epsilon$ is a half period.
\item $(E \times F_i)/G$, $G = {\Z}/4$ acts on $F_i$ by $x \mapsto ix$.
\item $(E \times F_i)/G$, $G = {\Z}/4 \oplus {\Z}/2$ acts on $F_i$ by
$x \mapsto ix$, $x \mapsto x + (1+i)/2$.
\item $(E \times F_{\rho})/G$, $G = {\Z}/3$ acts on $F_{\rho}$ by $x
\mapsto \rho x$.
\item $(E \times F_{\rho})/G$, $G = {\Z}/3 \oplus {\Z}/3$ acts on
$F_{\rho}$ by $x \mapsto \rho x$, $x \mapsto x + (1- \rho)/3$.
\item $(E \times F_{\rho})/G$, $G = {\Z}/6$ acts on $F_{\rho}$ by
$x \mapsto -\rho x$.
\end{enumerate}
\end{teo}

In fact, the characterization of hyperelliptic surfaces is one
of the two key steps of Enriques' classification of algebraic
surfaces, we have namely:

\begin{teo}
\label{class}  The complex surfaces $S$ with $K$ nef, $K^2 =0$, $p_g =0$,
and such
that either $S$ is algebraic with $q = 1$, or, more generally, such that
$b_1 = 2$, are hyperelliptic surfaces if and only if $kod(S) = 0$
(equivalently, iff
 the Albanese fibres are smooth of genus 1).
\end{teo}

I want now to illustrate the main technical ideas used for
the classification of real Hyperelliptic Surfaces.

\begin{df}
\label{hatg} The extended symmetry group $\hat{G}$ 
is the group generated by
$G$ and a lift 
$\tilde{\sigma}$ of $\sigma$.
\end{df}

Rerunning the classification theorem yields

\begin{teo}
\label{ISO1}  Let $(S, \sigma)$, $(\hat{S}, \hat{\sigma})$ be isomorphic 
real hyperelliptic surfaces : then the
respective extended symmetry groups $\hat{G}$ are the same and given 
 two Bagnera - De Franchis realizations
 $S = (E \times F)/G$, $\hat{S} = (\hat{E} \times
\hat{F})/G$, there is an isomorphism 
$\Psi: E \times F \rightarrow \hat{E} \times \hat{F}$, of product type, 
commuting with the action of $\hat{G}$, and inducing the given isomorphism
$\psi: S \stackrel{\cong} \rightarrow \hat{S}$.
Moreover, let $\tilde{\sigma} : E \times F \rightarrow E \times F$ be
a lift of $\sigma$. Then
the antiholomorphic map $\tilde{\sigma}$ is of product type.
\end{teo}

We need to give the list of all the possible groups $\hat{G}$.

\begin{lem}
\label{action}

Let us consider the extension
$$(*) \ 0 \rightarrow G \rightarrow \hat{G} \rightarrow {\Z}/2 = <\sigma>
\rightarrow 0.$$

We have the following possibilities for the action of $\sigma$ on $G$:
 in what follows the subgroup $T$ of $G$ will be the subgroup
 acting by translations on both factors.

\begin{enumerate}
\item If $G = {\Z}/q$, $ q= 2,4, 3, 6$ then $\sigma$ acts as $-Id$ on $G$
and
$\hat{G} = D_q$

\item If $G = {\Z}/2 \times {\Z}/2$, then either 

\begin{itemize}
\item (2.1)
$\sigma$ acts  as the identity
on $G$,  $(*)$ splits, and $\hat{G} = {\Z}/2 \times {\Z}/2 \times {\Z}/2$
or
\item (2.2) 
$\sigma$ acts  as the identity
on $G$,  $(*)$  does not split, and $\hat{G} = {\Z}/4 \times {\Z}/2$ (and
in this latter case the square of $\sigma$ is the generator of $T$) or
\item (2.3)
 $\sigma$ acts as $ \left( \begin{array}{cc} 1 & 0\\ 1 & 1
\end{array} \right)
$,   $(*)$ splits, $\hat{G} = D_4$, the dihedral group, and again the
 square of a generator of $\Z/4 \Z$ is the generator of $T$.
\end{itemize}
\item If $G = {\Z}/4 \times {\Z}/2$, then either
$\hat{G} = T \times D_4 \cong {\Z}/2 \times D_4$, or
$\hat{G} $ is isomorphic to the group $G_1:= <\sigma, g, t, \ | \sigma^2 =1, \ g^4 = 1, \ t^2 = 1, \ t \sigma = \sigma t, \ tg = gt, \ \sigma g = g^{-1} t \sigma >$, and its action on the second
elliptic curve $F$ is generated by the following  transformations:
$ \sigma (z) =  \bar{z} + 1/2$, $ g(z) = iz$, $t(z) = z + 1/2 (1+i)$.
[The group $G_1$ is classically denoted by $c_1$ (cf.Atlas of Groups)].

In particular, in both cases  $(*)$ splits.

\item If $G = {\Z}/3 \times {\Z}/3$, then we may choose a subgroup $G'$ of
$G$ such that
$\sigma$ acts as
$-Id \times Id$ on $ G' \times T = G$ and we have $\hat{G} = D_3 \times
{\Z}/3$.

\end{enumerate}
\end{lem} 

MAIN IDEAS IN THE CLASSIFICATION OF REAL HYPERELLIPTIC SURFACES : 

i) the first point is to show that the orbifold fundamental group
has a unique faithful
representation as a group of affine transformations with rational
coefficients

ii) second, to distinguish the several cases, we use byproducts of the
orbifold fundamental group, such as the structure of the
extended Bagnera de Franchis group $\hat{G}$, parities of the
(invariant $\nu$ of the) involutions
 in $\hat{G}$, their actions on the fixed point sets of
 transformations in $G$, the topology of the real part ,
 and, in some special cases, some more refined invariants
 such as the translation parts of all possible lifts of a given element in
$\hat{G}$ to the orbifold fundamental group.

 HOW TO CALCULATE THE TOPOLOGY OF THEIR REAL PART ?
 
We simply use the Albanese variety and what we said about the real parts
of elliptic curves, that they consist of $\nu \leq 2$ circles.
Since also the fibres are elliptic curves, we get up to $4$
circle bundles over circles, in particular the connected components are
either 2-Tori or Klein bottles. 

For every connected component $V$ of $Fix(\sigma) = S({\R})$
its inverse image
$\pi^{-1}(V)$ in $ E \times F$ splits as the $G$-orbit of any of its
connected components. Let $W$ be one such: then one can easily see that
there is a lift
$\tilde{\sigma}$ of
$\sigma$ such that $\tilde{\sigma}$ is an involution  and
$W$ is in the fixed locus of $\tilde{\sigma}$: moreover 
$\tilde{\sigma}$ is unique.

Thus the connected components of
$Fix(\sigma)$ correspond bijectively to the set $\mathcal{C}$ 
obtained as follows:
consider all the lifts $\tilde{\sigma}$ of $\sigma$ which are involutions 
and pick one representative $\tilde{\sigma_i}$ for each conjugacy class.

Then we let $\mathcal{C'}$ be the set of equivalence classes of connected
components
of
$\cup Fix(\tilde{\sigma_i})$, where two components $A$, $A'$ of
$Fix(\tilde{\sigma_i})$ are equivalent if and only if  there exists an
element $g
\in G$, such that $g(A) = A'$.

Let now $\tilde{\sigma}$ be an antiholomorphic involution which is a lift of
$\sigma$ and such that $Fix(\tilde{\sigma}) \neq \emptyset$.

Since $\tilde{\sigma}$ is of product type, $Fix(\tilde{\sigma})$ is a 
disjoint
union of $2^{a_1 + a_2}$ copies of ${\bf S}^1 \times {\bf S}^1$, where $a_i \in
\{0,1\}$.

In fact if an antiholomorphic involution $\hat{\sigma}$ on an elliptic curve
$C$ has fixed points, then $Fix(\hat{\sigma})$ is a disjoint union of $2^a$
copies of ${\bf S}^1$, where $a =1$ if the matrix of the action of
$\hat{\sigma}$ on $H_1(C, \Z)$ is diagonalizable, else $a =0$.

What said insofar describes the number of such components $V$;
in order to determine their nature, observe that $V = W/H$, where $W$ is as
before and $H \subset G$ is the subgroup such that $HW = W$.

Since the action of $G$ on the first curve $E$ is by translations, the
action on
the first ${\bf S}^1$ is always orientation preserving. Whence,
$V$ is a Klein bottle if and only if $H$ acts on the second ${\bf S}^1$ by
some
orientation reversing map, or, equivalently, iff $H$ has some fixed point on
the second
${\bf S}^1$.

In fact, we have  a restriction : let  $h \in H$ be a transformation having
a fixed point on the second ${\bf S}^1$. 
Since the direction of this ${\bf
S}^1$ is an eigenvector for the tangent
action of $h$, it follows that the tangent action is given by 
multiplication by
$-1$.

FINALLY , EXPLICIT TABLES SHOW THAT REALLY EVERYTHING CAN BE
NICELY WRITTEN DOWN AS IN THE BAGNERA  DE FRANCHIS LIST
(we refer to the original paper \cite{c-f2} for details).

\section{Tori and curve times a torus}

In this section we begin analysing in detail the case of the simplest
complex 
 manifolds which are $K(\pi, 1)$'s, namely the complex tori.
 
 It is well known that complex tori are parametrized by a connected
family (with smooth base space), inducing  all the small deformations (cf.
\cite{k-m}).

A first result we want to establish is the following 

\begin{teo}
Every deformation of a complex torus of dimension $n$ is a 
complex torus of dimension $n$.
\end{teo}

The proof  starts with a sequence of more or less standard arguments.
The first one is due to Kodaira ( \cite{koi}, Theorem 2, page 1392 of
collected works, vol. III)

\begin{lem}
On a compact complex manifold $X$ one has an injection

$$ H^0( d \hol_X ) \oplus  H^0( d  \bar{\hol_X })   
\rightarrow  H^1_{DR} (X,  \C).$$
\end{lem}

PROOF. Suffices to show the injection $ H^0( d \hol_X )   
\rightarrow  H^1_{DR} (X,  \R).$ by the map sending 
$$ \omega  \rightarrow \omega  + \bar{\omega } .$$
Else, there is a fucntion $f$ with $df = \omega  + \bar{\omega }$,
whence $ \partial f = \omega $ and therefore $ \bar{\partial} 
\partial f =
d (\omega) = 0$. Thus $f$ is pluriharmonic, hence constant by the maximum
principle. Follows that $\omega = 0$.
\hfill Q.E.D. \\

\begin{lem}
Assume that $\{X_t\}_{t \in \Delta}$ is a 1-parameter family of compact
complex manifolds over the 1-dimensional disk, such that there is a
sequence $ t_{\nu} \rightarrow 0$ with $X_{t_{\nu}}$ K\"ahler.

Then the  weak 1-K\"ahler  decomposition 
$$ H^1_{DR} (X_0,  \C) = H^0( d \hol_{X_0} ) \oplus  
H^0(\bar{d \hol_{X_0}} ),  $$
holds also on the central fibre $X_0$.
\end{lem}

PROOF. We have $ f: \Xi \rightarrow \Delta  $ which is proper and smooth,
and  $ f_* (\Omega^1_{\Xi |\Delta}) $  is torsion free, whence ($\Delta$ is
smooth of dimension $1$) it is locally free of rank $ q: = (1/2) \
b_1(X_0).$

In fact, there is ( cf. \cite{mum}, II 5) a complex of Vector Bundles on
$\Delta$,

$$  (*)\  \  E^0 \rightarrow E^1 \rightarrow E^2 \rightarrow ...E^n\ 
s.t. $$ 
$1) R^if_* (\Omega^1_{\Xi |\Delta}) $ is the i-th cohomology group of (*),
whereas 

$ 2)H^i (X_t, \Omega^1_{X_t}) $ is the i-th cohomology group of
$(*) \otimes \C_t .$

CLAIM. 

Thus there are holomorphic 1-forms $\omega_1(t), ...\omega_q(t)$
defined  in the inverse image $f^{-1} (U_0)$ of a neighbourhood 
$(U_0)$ of $0$ , and linearly independent for $ t \in U_0$.

PROOF OF THE CLAIM: assume that  $\omega_1(t), ...\omega_q(t)$ generate
the direct image sheaf  $ f_* (\Omega^1_{\Xi |\Delta}) $, but
$\omega_1(0), ...\omega_q(0)$ be linearly dependent. Then, w.l.o.g. we
may assume
$\omega_1(0)\equiv 0$, i.e. there is a maximal $m$ such that
$\hat{\omega}_1(t) : = \omega_1(t)/t^m$ is holomorphic. Then, since 
$\hat{\omega}_1(t) $ is a section of  $ f_* (\Omega^1_{\Xi |\Delta})
$, there are holomorphic functions $\alpha_i$ such that 
$\hat{\omega}_1(t)  = \Sigma_{i=1,..q} \alpha_i(t)\omega_1(t) $,
whence it follows that $\hat{\omega}_1(t) ( 1- t^m \alpha_1(t)) =
\Sigma_{i=2,..q} \alpha_i(t)\omega_1(t).$

This however contradicts the fact that $ f_* (\Omega^1_{\Xi |\Delta}) $ 
is locally free of rank $ q$.

REMARK:

Let $d_v$ be the vertical part of exterior differentiation, i.e., the
composition of $d$ with the projection $  (\Omega^2_{\Xi} \rightarrow
\Omega^2_{\Xi |\Delta}) $. By our asumption, $\omega_i(t_{\nu})$ is
$d_v$-closed, whence, by continuity,  also $\omega_i(0) \in H^0 (
d\hol_{X_{0}}).$ 

END OF THE PROOF:

It follows from the previous lemma that, being $b_1 = 2
q$, $ H^1_{DR} (X_0,  \C) = H^0( d \hol_{X_0}  ) \oplus  
H^0( d \bar{\hol}_{X_0} )$.

\hfill Q.E.D. \\ 

Recall that for a compact complex manifold $X$, the Albanese  Variety
$Alb(X)$ is the quotient of the complex dual vector space of $H^0( d \hol_X 
)$ by the minimal  closed complex Lie subgroup containing the image of
$H_1(X,
\Z)$.

The Albanese map
$\alpha_X : X \rightarrow Alb(X)$ is given as usual by fixing a
base point $x_0$, and defining $\alpha_X (x) $ as the integration on any 
path connecting $x_0$ with $x$.

One says the the Albanese Variety is {\bf good} if the image of $H_1(X, \Z)$
is discrete in $H^0( d \hol_X  )$, and {\bf very good } if it is a lattice.
Moreover, the {\em Albanese dimension of $X$ } is defined as the
dimension of the image of the Albanese map.

With this terminology, we can state an important consequence of our
assumptions.

\begin{cor}
Assume that $\{X_t\}_{t \in \Delta}$ is a 1-parameter family of compact
complex manifolds over the 1-dimensional disk, such that there is a
sequence $ t_{\nu} \rightarrow 0$ with $X_{t_{\nu}}$ K\"ahler, and
moreover such that $X_{t_{\nu}}$  has maximal Albanese dimension.

Then the central fibre $X_0$ has a very good Albanese Variety, and has
also maximal Albanese dimension.
\end{cor} 

PROOF. We use the fact (cf. \cite{cat4}) that, when the Albanese Variety is
good, then  the Albanese dimension of
$X$ is equal to $ max \{  i| \Lambda^i H^0( d \hol_X  ) \otimes 
\Lambda^i H^0( d \bar{\hol}_X  ) \rightarrow H^{2i}_{DR}
(X, \C)$ has non zero image $\}.$

If the weak 1-K\"ahler decomposition holds for $X$, then the Albanese
dimension of $X$ equals $(1/2) \   max \{  j| \Lambda^j H^1 (X, 
\C) $ has non zero image in   $ H^j (X, 
\C)\}.$  But this number is clearly invariant by
homeomorphisms.

Finally, the Albanese Variety for $X_0$ is very good since it is very
good for $X_{t_{\nu}}$ and the weak 1-K\"ahler decomposition holds for 
$X_0$.

\hfill Q.E.D. \\ 

\begin{oss}
As observed in (\cite{cat6}, 1.9), if a complex manifold $X$ has a
generically finite map to a K\"ahler manifold, then $X$ is bimeromorphic to
a  K\"ahler manifold. This applies in particular to the Albanese map.
\end{oss}

Recall that (cf. \cite{k-m}) a small deformation of a K\"ahler manifold
is again K\"ahler. This is however (cf. \cite{hir}) false for deformations
in the large, and shows that the following theorem ( which pretty much
follows the lines of corollary C of \cite{cat6}) is not entirely obvious
(it answers indeed in the affirmative a problem raised by Kodaira and
Spencer , page 464 of the paper \cite{k-s}, where the case
n=2 was solved in Theorem 20.2 )

\begin{teo}
A  deformation in the large of complex tori is a complex torus.

This follows form the following more precise statement: let $X_0$ be a
compact complex manifold such that its Kuranishi family of deformations
$\pi : \Xi
\rightarrow
\mathcal{B}$ enjoys the property that the set $\mathcal{B}(torus) : = \{b|
X_b $is isomorphic to a complex torus$\}$ has $0$ as a limit point.

Then $X_0$ is a complex torus.
\end{teo}

We will use the following "folklore"

\begin{lem}
Let $Y$ be a connected complex analytic space, and $Z$ an open set of $Y$
such that $Z$ is closed for holomorphic 1-parameter limits ( i.e., given
any holomorphic map of the 1-disk $f : \Delta \rightarrow Y$, if there is
a sequence $t_{\nu} \rightarrow 0$ with $f(t_{\nu} ) \in Z$, then also
$f(0) \in Z$). Then $Z=Y$.
\end{lem}

PROOF OF THE LEMMA. 

By choosing an appropriate stratification of $Y$ by smooth manifolds, it
suffices to show that the statement holds for $Y$ a connected manifold.

Since it suffices to show that $Z$ is closed,  let $P$ be a
point in the closure of $Z$, and let us take coordinates such that a
neighbourhood of $P$ corresponds to a compact polycylinder $H$ in $\C^n$.

Given a point in $Z$, let $H'$ be a maximal coordinate polycylinder
contained in $Z$. We claim that $H'$ must contain $H$, else, by the
holomorphic 1-parameter limit property, the boundary of $H'$ is contained in
$Z'$, and since $Z$ is open, by compactness we find a bigger polycylinder
contained in $Z$, a contradiction which proves the claim.

\hfill Q.E.D. for the Lemma.\\ 

PROOF OF THE THEOREM.

It suffices to consider a 1-parameter family ($\mathcal{B} = \Delta$)
whence we may assume w.l.o.g. that the weak 1-Lefschetz decomposition
holds for each
$t \in \Delta$.

By integration of the holomorphic 1-forms on the fibres ( which are
closed for $t_{\nu}$ and for $0$), we get a family of Albanese maps
$\alpha_t : X_t
\rightarrow J_t$, which fit together in a relative map $a : \Xi
\rightarrow J$ over $\Delta$ ( $J_t $ is the complex torus
$(H^0(d \hol_{X_t}))^{\vee}/ H_1(X_t, \Z)$).

Apply once more  vertical exterior differentiation to the forms 
$\omega_i(t)$: $ d_v (\omega_i(t)$  vanishes identically on $X_0$ and on
$X_{t_{\nu}}$, whence it vanishes identically in a neighbourhood of
$X_0$, and therefore these forms $\omega_i(t)$ are closed for each $t$.

Therefore our map $a : \Xi \rightarrow J$ is defined everywhere and it is an
isomorphism for $t = t_{\nu}.$

Whence, for each $t$, $\alpha_t$ is surjective and has degree 1.

To show that $\alpha_t$ is an isomorphism for each $t$ it suffices
therefore to show that $a$ is finite. 

Assume the contrary: then there is a ramification divisor $R$ of $a$,
which is exceptional (i.e., if $B = a(R)$, then $dim B < dim R$).

By our hypothesis $\alpha_t$ is an
isomorphism for $t = t_{\nu}$, thus $R$ is contained in a union of fibres,
and since it has the same dimension, it is a finite union of fibres.
But if $R$ is not empty, we reach a contradiction, since then there are
some $t$'s such  that $\alpha_t$ is not surjective.

\hfill Q.E.D. \\ 

\begin{prop}
Assume that $X$ has the same integral cohomology algebra of a complex
torus and that   $ H^0( d \hol_X  )$ has dimension equal to $n =dim(X)$.
Then $X$ is a complex torus
\end{prop}

PROOF.

Since $b_1(X) = 2n$, it follows from 4.2 that the weak 1-Lefschetz
decomposition holds for $X$ and that the Albanese Variety of $X$ is very
good.

That is, we have the Albanese map $\alpha_X : X \rightarrow J$, where 
$J$ is the complex torus $J= Alb(X)$. We want to show that the Albanese
map is an isomorphism. It is a morphism of degree $1$, since $\alpha_X $
induces an isomorphism between the respective fundamental classes of
$H^{2n}(X,\Z) \cong H^{2n}(J,\Z)$.

There remains to show that the bimeromorphic morphism $\alpha_X $ is finite.

To this purpose, let $R$ be the ramification divisor of $\alpha_X : X
\rightarrow J$, and $B$ its branch locus, which has codimension at least
$2$. By means of blowing ups of $J$ with non singular centres we can
dominate
$X$ by a K\"ahler manifold $g: Z \rightarrow X $ (cf. \cite{cat6} 1.8, 1.9).

Let $W$ be a fibre of $\alpha_X$ of positive dimension such that
$g^{-1} W$ is isomorphic to $W$. Since $Z$ is K\"ahler, $g^{-1} W$ is not
homologically trivial, whence we find a differentiable submanifold $Y$ of
complementary dimension which has a positive intersection number with it.
But then, by the projection formula , the image $g_* Y$ has positive
intersection with
$W$, whence $W$ is also not homologically trivial. However, the image of
the class of $W$ is $0$ on $J$, contradicting that $\alpha_X$ induces an
isomorphism of cohomology ( whence also homology) groups.

\hfill Q.E.D. \\ 

\begin{oss}
The second condition holds true as soon as the complex dimension $n$ is
at most $2$. For $n=1$ this is well known, for $n=2$ this is also known,
and due to Kodaira (\cite{koi}): in fact for $n=2$ the holomorphic 1-forms
are closed, and moreover  $ h^0( d \hol_X  )$ is at least $[ (1/2)
b_1(X)]$.

For $n\geq 3$, the real dimension of $X$ is greater than $5$, whence, by
the s-cobordism theorem (\cite{maz}), the assumption that $X$ be
homeomorphic to a complex torus is equivalent to the assumption that  $X$ be
diffeomorphic to a complex torus. 

Andr\'e Blanchard (\cite{bla1}) constructed in the early 50's an example of
a non K\"ahler complex structures on the product of a rational curve with a
two dimensional complex torus. In particular his construction (cf.
\cite{somm}), was rediscovered by Sommese, with a more clear and
 more general presentation, who pointed out that in this way one would
produce exotic complex structures on complex tori. 

The Sommese-Blanchard examples (cf. \cite{ue}, where we learnt about them)
are particularly relevant to show that the plurigenera of non K\"ahler
manifolds are not invariant by deformation. In the following  section 5 we
will adopt the presentation of 
\cite{cat6} ) : the possible psychological reason why we had forgotten
about these 3-folds (and we are very thankful to Andrew Sommese
for pointing out their relevance) is that these complex structures do not
have a trivial canonical bundle. Thus remains open the following

 QUESTION: let $X$ be a  compact complex manifold of dimension $n \geq
3$  and with trivial canonical bundle such that $X$ is diffeomorphic to a
complex torus : is then $X$ a complex torus ?
\end{oss}

The main problem is to show the existence of holomorphic 1-forms, so
it may well happen that also this question has a negative answer.

On the other hand, as was already pointed out, the 
examples  of Blanchard, Calabi and Sommese show  that the answer to a
similar question is negative also in the case of a product $C \times T$,
where
$C$ is a curve of genus $g\geq 2$, and $T$ is a complex torus.

For this class of varieties we still have the same result as for tori,
concerning global deformations.

\begin{teo}
A  deformation in the large of products $C \times T$ of a curve of
genus $g\geq 2$ with a complex torus is again  a product of this type.

This clearly follows from the more precise statement: let  $g\geq 2$ a fixed
integer, and assume that
$X_0$ be a compact complex manifold such that its Kuranishi family of
deformations
$\pi : \Xi \rightarrow \mathcal{B}$ enjoys the property that the set
$\mathcal{B}" : = \{b| X_b$ is isomorphic to the product of a curve of
genus $g$ with a complex torus $\}$ has
$0$ as a limit point.

Then $X_0$ is isomorphic to a product $C \times T$ of a curve of
genus $g$ with a complex torus $T$.
\end{teo}

PROOF.
Observe first of all that $\mathcal{B}"$ is open in the Kuranishi family
$\mathcal{B}$,
by the property that the Kuranishi family induces a versal family in each
neighbouring point, and that every small deformation of a product
$C\times T$ is of the same type.

Whence, by lemma  $6.7$ we can limit ourselves to consider the situation
where $\mathcal{B}$ is a 1-dimensional disk, and every fibre $X_t$ is,
for $t \neq 0$, a product of the desired form.

STEP I. There is a morphism $F: \Xi \rightarrow \mathcal{C}$, where 
$\mathcal{C} \rightarrow \mathcal{B}$ is a smooth family of curves of genus
$g$.

{\bf Proof of step I.} 

We use for this purpose the isotropic subspace theorem
of \cite{cat4}, observing that the validity of this theorem does not require
the full hypothesis that a variety $X$ be K\"ahler, but that 
weaker hypotheses, e.g. the  weak 1-K\"ahler and 2-K\"ahler decompositions
do indeed suffice.

The first important property to this purpose is that the cohomology algebra
$H^{*}(X_t,\C)$ is generated by $H^1(X_t,\C)$.

The second important property is that for a product $C\times T$
as above, the subspace $p_1^*(H^1(C,\C))  $ is the unique maximal subspace
$V$, of dimension $2g$, such that the image of $\Lambda^3(V) \rightarrow
H^3(C \times T,\C)  $ is zero.

This is the algebraic counterpart of the geometrical fact that the first
projection is the only surjective morphism with connected fibre of
$C\times T$ onto a curve $C'$ of genus $\geq2$.

From the differentiable triviality of our family $\Xi \rightarrow \Delta$
follows that we have a uniquely determined such subspace $V$ of the
cohomology of $X_0$, that we may freely identify to the one of each $X_t$.

Now, for $t\neq 0$, we have a decomposition $V = U_t \oplus \bar{U}_t$,
where $ U_t = p_1^*(H^0(\Omega^1_{C_t}))  .$

By compactness of the Grassmann variety and by the  weak 1-K\"ahler
decomposition in the limit, the above decomposition also holds for $X_0$,
and $U_0$ is a maximal isotropic subspace in  $ H^0( d \hol_{X_0} )$.
The Castelnuovo de Franchis theorem applies, and we get the desired
morphism $F: \Xi \rightarrow \mathcal{C}$ to a family of curves.

STEP II. 
We  produce now a morphism $G: \Xi \rightarrow \mathcal{T}$
where $\mathcal{T} \rightarrow \mathcal{B}$ is a family of tori.

Proof of step II.

Let $\alpha : \Xi \rightarrow \mathcal{J}$ be the family of Albanese maps.
By construction, or by the universal property of the Albanese map, we
have that $F$ is obtained by taking  projections of the Albanese maps,
whence we have a factorization $\alpha : \Xi \rightarrow \mathcal{J}
\rightarrow \mathcal{C} \rightarrow \mathcal{J"}$, where $\mathcal{J"}
\rightarrow \mathcal{B}$ is the family of Jacobians of $\mathcal{C}.$

The desired family of tori $\mathcal{T} \rightarrow \mathcal{B}$ is
therefore the family of kernels of $\mathcal{J} \rightarrow \mathcal{J"}.$

To show the existence of the morphism $G$, observe that for $t\neq 0$, an
isomorphism of $X_t$ with $C_t \times T_t$ is given  by a projection
$ G_t: X_t \rightarrow T_t$ which, in turn, is provided (through
integration) via a complex subspace
$W_t$ of
$H^0(\Omega^1_{X_t})$ whose real span $H_t =W_t \oplus \bar{W_t} $ is
 $\R$-generated by a subgroup of
$H^1(X_t,\Z)$.

Although there are several choices for such a subgroup, we have at most a
countable choice of those. Since for each $t$ there is such a
choice, by Baire's theorem there is a choice which holds for each $t\neq
0$. Let us make such a choice of $H = H_t$ $\forall t$: then the
corresponding subspace
$W_t$ has a limit in $H$ for
$t=0$, and since the weak 1-Lefschetz decomposition holds for $X_0$, and
this limit is a direct summand for $U_t$, we easily see that this limit is
unique, and  the desired morphism $G$ is therefore obtained.

STEP III. The fibre product $F\times_{\mathcal{B}} G$ yields an
isomorphism onto ${\mathcal{C}} \times_{\mathcal{B}} {\mathcal{T}}$
for each fibre over $t\neq 0$, but since also the fibre over $0$ is a
torus fibration over $C_0$, and the cohomology algebras are all the same,
it follows that we have an isomorphism also for $t=0$.

\hfill Q.E.D. \\

We end this section by  giving an application of theorem 4.6.

\begin{teo}
Let $(X, \sigma)$ be a real variety which is a deformation of real tori:
then also $(X, \sigma)$ is a real torus.

The orbifold fundamental groupo sequence completely determines the
differentiable type of the pair $(X, \sigma)$.

The real
deformation type is also  completely determined by the orbifold
fundamental group of $(X, \sigma)$.

In turn, the orbifold fundamental group, if $s$ is the linear integral
transformation of $H_1(X, \Z)$ obtained by conjugation with $\sigma$, is
uniquely determined by the integer $r$ which is the rank of the matrix
$(s- Id )$$(mod 2)$ acting on $H_1(X, \Z/2)$ and, if $r\neq n=dim X$, by
the property whether the orbifold exact sequence does or does not split (if
$r=n$ it always splits).
\end{teo}

PROOF. 

The first statement is a direct consequence of theorem 4.6. 

Concerning the second statement, let  $\tilde{\sigma}$ be  a lifting of the
antiholomorphic involution $\sigma$ on the universal cover $\C^n$: then
it follows as usual that the first derivatives are bounded, whence
constant, and $\tilde{\sigma}$ is an affine map. We may moreover
assume  $\tilde{\sigma}$ to be linear in the case where $\sigma$ admits a
fixed point (in this case the orbifold exact sequence splits).

Let $\Lambda$ be the lattice $H_1(X, \Z)$: so we have first of all that
$\tilde{\sigma}$ is represented by the affine map of $\Lambda \otimes \R$,
$x \rightarrow s(x) + b$, where $s$ is the isomorphism of $\Lambda$ to
itself, given by conjugation with $\sigma$.

Thus $s$ is an integral matrix whose square is the identity, and  it is
well known ( cf. e.g. [Cat 9], lemma 3.11) that we can thus split 
$\Lambda = U \oplus V \oplus W^+  \oplus W^-$, where $s$ acts by
$ s(u,v,w^+,w^-) = (v,u,w^+,- w^-).$

Since the number of $(+1)$ -eigenvalues is equal to the number of
$(-1)$ -eigenvalues, $\tilde{\sigma}$ being antiholomorphic, we get that 
$W^+  , W^-$ have the same dimension $n-r$ (whence $dim U= dim V = r$).

Since the square of $\sigma$ is the identity, it follows that 
$\tilde{\sigma}^2$  is an integral translation. This condition boils down
to $ s b + b \in \Lambda$. If we write $b= (b_1,b_2,b^+,b^-)$, the
previous condition means that $b_1 +b_2, 2b^+$ are integral vectors.

In order to obtain a simple normal form, we are allowed to take a
different lift $\tilde{\sigma}$, i.e., to consider an affine map of the
form $x \rightarrow s(x) + b + \lambda$, where $\lambda  \in \Lambda$,
and then to take another point $c \in \Lambda \otimes \R$ as the origin.

Then we get an affine transformation of the form $y \rightarrow sy + b
+ \lambda + s c - c$. The translation vector is thus 
$ (b_1,b_2,b^+,b^-) + (\lambda_1,\lambda_2,\lambda^+,\lambda^-) +
(c_1 - c_2,c_2 - c_1, 0, 2 c^-)$. 

We choose $  c^- = -(1/2) b^-$ and $\lambda^- = 0$, $ c_2 - c_1 = b_1
+\lambda_1$, $- \lambda_2 = b_1 +b_2 - \lambda_1$, and then we choose for
$\lambda^+$ the opposite of the integral part of $b^+$.

We obtain a new affine map $y \rightarrow sy + b'$ where $b' =
(0,0,{b'}^+,0)$ and all the coordinates of ${b'}^+$ are either $0$ or
$1/2$. So, either $b' =0$, or we can assume that $2 {b'}^+$ is the first
basis element of $W^+$.

We obtain thus two cases :

\begin{itemize}
\item
1) There exists a lift $\tilde{\sigma}$ represented by the linear map $s$
for a suitable choice of the origin.
\item
2) There exists $\tilde{\sigma}$ represented, for a suitable choice of the
origin, and for the choice of a suitable basis of $\Lambda$, by the affine
map $y \rightarrow sy + 1/2 e_1^+$ where $e_1^+$ is the first basis
element of $W^+$.

\end{itemize}

REMARK: Case 1) holds if and only if there is a fixed point of $\sigma$.
Whereas, in case 2) we can never obtain that ${b'}^+$ be zero, therefore
the square of any lift $\tilde{\sigma}$ is always given by a
translation with  third coordinate $2 {b'}^+ \neq 0$, so the orbifold exact
sequence does not split.

We conclude that 1) holds if and only if the orbifold exact
sequence  splits.

By the  normal forms obtained above, it follows that the integer $r$ and
the splitting or non splitting property of the orbifold fundamental group
sequence not only determine completely the orbifold fundamental group
sequence, but also its  affine representation.

To finish the proof of the theorem, assume that the normal form of
$\tilde{\sigma}$ is preserved, whence also the affine representation of the
orbifold fundamental group.

We are now looking for al the translation invariant complex structures
which make the transformation $\tilde{\sigma}$ antiholomorphic.
As before in the case of the elliptic curves, we look for the $2n
\times 2n$ matrices
$J$ whose square $J^2 = - Id$, and such that $Js = -s J$. 
The latter condition implies that $J$ exchanges the eigenspaces of $s$ in 
$\Lambda \otimes \R.$

We can calculate the matrices $J$ in a suitable $\R$-basis of $\Lambda
\otimes \R$ where $s$ is diagonal, i.e., $s(y_1, y_2) =(y_1,- y_2)$.
Then $J(y_1, y_2) = (A  y_2, B  y_1 )$ and $Js = -s J$ is then satisfied.
The further condition $J^2 = - Id$ is equivalent to $AB=BA= -Id$, i.e.,
to $ B = - A^{-1}$.
Therefore, these complex structures are parametrized by $GL(n, \R)$,
which indeed has two connected components, distinguished by the sign of
the determinant.

Since however we already observed that $\sigma$ provides an isomorphism
between $(X, \sigma)$ and $(\bar{X}, \sigma)$,  and $\bar{X}$ corresponds
to the complex structure $-J$, if $J$ is the complex structure for $X$,
we immediately obtain that $A$ and $-A$ give isomorphic real varieties.

We are set for $n$ odd, since $det (-A) = (-1)^n det A$. 

In the case where $n$ is even, observe that two matrices $A$ and $A'$
yield isomorphic real tori if and only if there is a matrix $D \in
GL(\Lambda)$ such that $D$ commutes with $s$, 
  and such that $D$ conjugates $J$ to $J'$ (then the diffeomorphism is
orientation preserving for the orientations respectively induced by the 
complex structures associated to $A$, resp. $A'$). 
  
  Thus $D$ respects the
eigenspaces of $s$, whence $D(y_1, y_2) =(D_1 y_1,D_2 y_2)$, and $A$ is
transformed to $ A' = D_2 A (D_1)^{-1}$. Whence, we see that the 
sign of the determinant of
$A'$ equals the one of $ (detA)  (det D ) $, whence we can change  sign
in any case by simply choosing $D$ with $det D = -1$.

\hfill Q.E.D. \\ 

\section{The Blanchard -Calabi threefolds}

The  Sommese-Blanchard examples (\cite{bla1}, \cite{bla2}, \cite{somm})
provide non K\"ahler complex structures
$X$ on manifolds diffeomorphic to a product $C \times T$, where $C$ is a
compact complex curve and $T$ is a 2-dimensional complex torus.

In fact, in these examples, the projection $X \rightarrow C$ is holomorphic 
and all the fibres are 2-dimensional complex tori.

Also Calabi (\cite{cal})   showed that there are  complex
structures  on a product $C\times T$ (that all these structures are 
non K\"ahlerian  follows also by the arguments of the previous theorem, 
else  they would produce a complex product structure
$C' \times T'$).

The result of Calabi is the following

\begin{teo} (CALABI)
Let $C$ be a hyperelliptic curve of  odd genus $g\geq 3$, and let $T$ be
a two dimensional complex torus. Then the differentiable manifold $C
\times T$ admits a complex structure with trivial canonical bundle.
\end{teo}

We shall try to show that the construction of Calabi, although formulated
in a different and very interesting general context, yields indeed a very
special case of the construction of Sommese-Blanchard,
which may
instead be formulated and generalized as follows

\begin{teo} (BLANCHARD-CALABI Jacobian 3-FOLDS)

Let $C$ be curve of   genus $g\geq 0$, and let $W$ be a rank $2$
holomorphic vector bundle admitting four holomorphic sections
$\sigma_1,\sigma_2,\sigma_3,\sigma_4$ which are everywhere
$\R$-linearly independent (for instance, if $W = L \oplus L$, where
$H^0(L)$ has no base points, then four sections as above do exist).

Then the quotient $X$ of $W$ by the
$\Z^4$-action acting fibrewise by translations : $ w \rightarrow  w +
\Sigma_{i=1,..4}\ n_i \sigma_i$ is a complex manifold diffeomorphic to 
the differentiable manifold $C \times T$,  where $T$ is a two dimensional
complex torus, and will be called a {\bf Jacobian Blanchard-Calabi 3-fold}.

The canonical divisor of $X$ equals $
K_X =
\pi^* (K_C - det W) $, where $\pi : X \rightarrow C$ is the canonical
projection. Moreover, $X$ is K\"ahler if and only if the bundle $W$ is
trivial (i.e., iff $X$ is  a holomorphic product $C \times T$).

Indeed one
has $h^0 (\Omega^1_X) = g $  unless $W$ is trivial, while, if the vector
bundle $V$ is defined through the exact sequence 
$$ 0 \rightarrow V \rightarrow (\hol_C)^4
\rightarrow  W \rightarrow 0,$$
then $h^1 (\hol_X) = g + h^0 (V^{\vee}) $.
\end{teo}

PROOF.
The four holomorphic sections $\sigma_1,\sigma_2,\sigma_3,\sigma_4$ make
$W$ a trivial $\R$-vector bundle, and with this trivialization we obtain
that $X$ is diffeomorphic to the product of $C$ with a real four
dimensional torus.

Let us now show that any vector bundle $W = L \oplus L$, where
$H^0(L)$ has no base points, admits such sections
$\sigma_1,\sigma_2,\sigma_3,\sigma_4$.

In fact, we have  $ s_1, s_2 \in H^0(L) $ without common zeros, so $s:=
^t(s_1, s_2)$ is a nowhere vanishing section of $W = L \oplus L$.

Use now the fact that $GL(2, \C)$ operates on $ L \oplus L$, and that,
identifying $\C^2$ with the field $H$  of Hamilton's quaternions, 
then $GL(2, \C)$ contains the finite quaternion group $\mathcal {H} =
\{+1, -1,  +i, - i, +j, - j, +ij, - ij \}$, thus it suffices to 
define $\sigma_1:= s, \sigma_2:= is, \sigma_3:= js, \sigma_4:= ijs.$

In other more concrete words $\sigma_2 :=
^t(is_1, -is_2), \sigma_3:= ^t(-s_2, s_1),\sigma_4 :=
^t(is_2, is_1)$.

We have the Koszul complex associated to $(s_1, s_2)$: 

$ 0 \rightarrow L^{-1} \rightarrow (\hol_C)^2 \rightarrow  L \rightarrow 0.$

Now, an easy calculation shows that, defining $V$ as the kernel subbundle
of the linear map given by $\sigma:= (\sigma_1,\sigma_2,\sigma_3,\sigma_4)$

$ 0 \rightarrow V \rightarrow (\hol_C)^4 \rightarrow  W \rightarrow 0,$

then in the special case above we have $V = L^{-1} \oplus L^{-1}$. 

Returning to the general situation,
since now on a torus $T= W'/\Gamma$ there are canonical isomorphisms of
$H^0(\Omega^1_T)$ with the space of linear forms on the vector space $W'$,
and of $H^1(\hol_T)$ with the quotient vector space Hom $(\Gamma, \C )$/
$H^0(\Omega^1_T)$, we obtain immediately 

1) the exact sequence 

$$ 0 \rightarrow \pi^*(\Omega^1_C)  \rightarrow \Omega^1_X \rightarrow
(W^{\vee}) \rightarrow 0   $$

2) an isomorphism $ \mathcal{R}^1 \pi_* (\hol_X) \cong (V^{\vee}).$

From the Leray spectral sequence follows immediately that 
 $h^1 (\hol_X) = g + h^0 (V^{\vee}) $.
 
 If $X$ is K\"ahler, then, since the first Betti number of $X$ equals
$4+2g$, then it must hold that $h^0 (\Omega^1_X) = g + 2$, in particular
$h^0 (W^{\vee}) \geq 2$ .

But $(W^{\vee})$ is a subbundle of a trivial bundle of rank $4$,
whence two linearly independent sections of $(W^{\vee})$ yield 
a composition $(\hol_C)^2 \rightarrow (W^{\vee}) \rightarrow
(\hol_C)^4   $ whose image is a trivial bundle of rank $1$ or $2$.

But in the former case the two sections would not be $\C$-linearly
independent, whence the image must be a trivial bundle of rank  $2$ and
taking determinants of the composition, we get that $(W^{\vee})$ is trivial
if  $h^0 (W^{\vee}) \geq 2$ .

Similarly, assume that $h^0 (W^{\vee}) =1$ : then $(W^{\vee})$ has a
trivial summand $I$ of rank $1$ . Then  we have a direct sum $(W^{\vee}) = I
\oplus Q$  and dually $W$ is a direct sum $\hol_C \oplus  Q^{-1}$. 
We use now the hypothesis that there are four holomorphic sections  of
$W$ which are $\R$-independent at any point: it follows then that $Q^{-1}$
admits two holomorphic sections  which are everywhere $\R$-independent.
Whence, it follows that also $Q^{-1}, Q$ are trivial.

\hfill Q.E.D. \\

\begin{df}
Given a Jacobian Blanchard-Calabi 3-fold  $ \pi : X \rightarrow C$,
any $X$-principal homogeneous space $\pi' :Y \rightarrow C$ will be called a
Blanchard-Calabi 3-fold.
\end{df}

\begin{cor}
The space of complex structures on the product of a curve $C$ with a four
dimensional real torus has unbounded dimension. 
\end{cor}

PROOF.

Let $G$ be the  four dimensional Grassmann variety $G(1,3)$: observe that
the datum of an exact sequence 
$$ 0 \rightarrow V \rightarrow (\hol_C)^4 \rightarrow  W \rightarrow 0$$
is equivalent to the datum of a holomorphic map $f: C\rightarrow G$,
since for any such $f$ we let $V,W$ be the respective pull backs of the
universal subbundle $U$ and  of the quotient bundle $Q$ (of course, for the
Blanchard-Calabi
 construction one needs the further open condition that the four
induced sections of $W$ be $\R$-linearly independent at each point).

We make throughout the assumption that  $f$ is not constant, i.e.,  $W$
is not trivial. Then we have the following exact sequence:
$$ 0 \rightarrow \Theta_C \rightarrow (f)^* \Theta_G  \rightarrow  N_f
\rightarrow 0 ,$$
where $N_f$, the normal sheaf  of the morphism $f$, governs the deformation
theory of the morphism $f$, in the sense that the tangent space to Def(f)
is the space $H^0(N_f)$, while the obstructions lie in $H^1(N_f)$.

By virtue of the fact that $\Theta_G  = Hom(U, Q)$, and of the
cohomology sequence associated to the aboce exact sequence, we get 
$$ 0 \rightarrow H^0(\Theta_C ) \rightarrow H^0(V^{\vee} \otimes W) 
\rightarrow  H^0(N_f) \rightarrow H^1(\Theta_C ) \rightarrow H^1(V^{\vee}
\otimes W) \rightarrow  H^1(N_f) 
\rightarrow 0 ,$$ and we conclude that the deformations of the map are
unobstructed provided $H^1(V^{\vee} \otimes W) = 0$.

This holds, in the special case where $W = L \oplus L$, if the degree $d$
of $L$  satisfies $d \geq g$, since then 
$H^1(V^{\vee} \otimes W) = H^1((L \oplus L) \otimes (L \oplus L) ) = 0$.

If $d \geq d$ the dimension of the space of deformations of the map $f$
is given by $ 3g - 3 + 4 h^0(2L) = 4 d + 1-g $, and this number clearly
tends to infinity together with $ d = deg(L)$.

On the other hand, we want to show that the deformations of the map $f$
yield effective deformations of $X$ as a Lie group principal fibration. This
can be seen as follows.

Consider the exact sequence 
$$ 0 \rightarrow \pi^* (W) \rightarrow  \Theta_X   \rightarrow  \pi^*
(\Theta_C)  \rightarrow 0 $$
and the derived direct image sequence  
$$ 0 \rightarrow (W) \rightarrow  \pi_* \Theta_X   \rightarrow 
(\Theta_C)  \rightarrow $$
$$  \rightarrow (V^{\vee} \otimes W) \rightarrow  \mathcal{R}^1 \pi_*
\Theta_X  
\rightarrow  (\Theta_C)\otimes V^{\vee}  \rightarrow $$
$$ \rightarrow (det(W) \otimes W ) \rightarrow  \mathcal{R}^2 \pi_*
\Theta_X  
\rightarrow  (\Theta_C)\otimes det(W ) \rightarrow 0,$$
where we used that  $\mathcal{R}^2 \pi_*(\hol_X) \cong \Lambda^2(V^{\vee})
\cong det(W).$

Notice that $ \pi_* \Theta_X $ is a vector bundle on the curve $C$, of rank
either $2$ or $3$.
In the latter case, since its image $M$ in $(\Theta_C)$ is saturated 
(i.e., $(\Theta_C)/M$ is torsion free), then $M = (\Theta_C)$ and all the
fibres of $\pi$ are then biholomorphic.

In this case, for any tangent vector field on $C$,
which we identify to the $0$-section of $W$, we  get a corresponding
tangent vector field on $W$ by  using the fiberwise simply transitive
action of $\R^4$ corresponding to the four chosen holomorphic sections.
Since all the fibres are isomorphic, the corresponding tangent vector field
on $X$ is holomorphic, and the bundle $W$ is then trivial.

Since we are assuming $W$ not to be trivial, we get thus an exact sequence
$$ 0   \rightarrow (\Theta_C)  \rightarrow 
 (V^{\vee} \otimes W) \rightarrow  \mathcal{R}^1 \pi_*
\Theta_X  $$
By the Leray spectral sequence follows that we have that
$H^0(V^{\vee} \otimes W) / H^0((\Theta_C))$ injects into 
$H^0(\mathcal{R}^1 \pi_*
\Theta_X  )$ which is a direct summand of $H^1(\Theta_X )= 
H^0(\mathcal{R}^1 \pi_*
\Theta_X  ) \oplus H^1(W).$

Therefore the smooth space $Def(f)$ of deformations of the map $f$ embeds
into the Kuranishi space $Def(X)$ of deformations of $X$.

The space $H^1(W)$  is the classifying space for principal $X$
homogeneous spaces  ( cf. \cite{ko-2-3}, 
\cite{shaf}) and therefore we see
that this subspace of $H^1(\Theta_X )$ corresponds to actual
deformations.

However for $W$ of large degree we get $H^1(W)=0$, so then the
Blanchard-Calabi 3-folds coincide with the Jacobian Blanchard-Calabi
3-folds.

\hfill Q.E.D. \\

\begin{oss}
1) In the Blanchard-Calabi examples one gets a trivial canonical bundle 
iff $ det(W) \equiv K_C$. This occurs in the special case where 
$W = L \oplus L$ and $L$ is a thetacharacteristic such that
$H^0(L)$ is base point free. In particular, we have the  Calabi
examples where $C$ is hyperelliptic of odd genus $g$ and $L$ is the
$1/2 (g-1)^{th}$ power of the hyperelliptic line bundle of $C$.

2) Start with a Sommese-Blanchard 3-fold with trivial canonical
bundle and deform the curve $C$ and the line bundle $L$ in such a way that
the canonical bundle becomes the pull back of a non torsion element of
$Pic(C)$: then this is the famous example that the Kodaira dimension is not
deformation equivalent for non K\"ahler manifolds ( cf. \cite{ue}).

3) In the previous theorem we have followed  the approach of \cite{cat6},
correcting however a wrong formula for $\mathcal{R}^1 \pi_* \hol_X$ (which
would give the dual vector bundle).

4) The approach of Sommese is also quite similar, and related to quaternion
multiplication, which is viewed as the $\C$-linear map $\psi : H = \C^2
\rightarrow Hom_{\R} ( H = \R^4 , H = \C^2)$, defined by 
$\psi (q) (q') = q q'$.

Whence, quaternion multiplication provides four sections of
$\hol_{\PP^1}(1)^2$ which are $\R$-linearly independent at each point.

Sommese  explores this particular situation as an example of the
natural fibration occurring in the more general context of the
so called quaternionic manifolds (later on, the word "twistor fibration",
for the other projection, has become more fashionable).

Sommese moreover observes that given a line bundle $L$ such that $H^0(L)$ is
base point free, the choice of two independent sections yields a
holomorphic map to
$\PP^1$ and   one can pull back $\hol_{\PP^1}(1)^2$ and the four
sections:  obtaining exactly the same situation we described above.

We learnt from Sommese that Blanchard knew  these examples too (they are
also described in  \cite{ue}). 
\end{oss}

We will indeed prove  in the sequel a much more precise statement
concerning  the deformations of Blanchard-Calabi 3-folds: we need for this
purpose the following 

\begin{df}
A Blanchard-Calabi 3-fold is said to be developable if the corresponding
ruled surface is developable, or, in other words , if and only if the
derivative of the corresponding map $f: C\rightarrow G$ yields at each point
$p$ of $C$ an element of
$(V^{\vee} \otimes W)_p$ of rank $\leq 1$. 
\end{df}
\begin{oss}
Observe in fact that (cf. \cite{a-c-g-h} C-9, page 38) since $f$
is non constant, and if there is no point of $\PP^3$ contained in  each
line of the corresponding family of  lines, $f$ is the associated map to a
holomorphic map $F :  C \rightarrow \PP^3$, so the union of the family 
of lines is the tangential developable of the image curve of the mapping
$F$.

If such a point exists, it means that there is an effective divisor $D$
on the curve $C$, and the bundle
$W$ is given as an extension 
$$ 0 \rightarrow (\hol_C)(-D) \rightarrow (\hol_C)^3 \rightarrow  W
\rightarrow 0$$
and moreover a fourth section of $W$ is given (which together with the
previous three provides the desired trivialization of the underlying  real
bundle).
\end{oss}

\begin{oss}
By a small variation of a theorem of E. Horikawa concerning the
deformations of holomorphic maps, namely Theorem 4.9 of
\cite{cat4}, we obtain in particular that, under the assumption $H^0
((\Theta_C)\otimes V^{\vee} ) =0$,  we have a smooth morphism $Def (\pi)
\rightarrow Def(X)$ . This assumption is however not satisfied in the case
of Sommese-Blanchard 3-folds, at least in the case where the degree of $L$
is large.
\end{oss}

\begin{teo}
Any small deformation of a non developable Blanchard-Calabi 3-fold such
that $H^1(V^{\vee} \otimes W) = 0$ is again a Blanchard-Calabi 3-fold. 
In particular, a small deformation of a 
Sommese-Blanchard 3-fold with $L$ of degree $d \geq
g$ is a Blanchard-Calabi 3-fold. 
\end{teo}

PROOF.

In this context recall once more (cf. \cite{ko-2-3})
that the infinitesimal deformations corresponding to $ H^1(W)$ correspond
to the deformations of $X$ as a principal homogenous space fibration over
the principal Lie group fibration $X$, and these are clearly unobstructed.
Indeed, for the Sommese-Blanchard examples, it will also hold
$H^1(W)=0$, if the degree $d$ of $L$ is large enough.

Assume that the vanishing  $H^1(V^{\vee} \otimes W) = 0$ holds:
then the deformations of the map $f$ are unobstructed, and
it is then clear that we obtain a larger family of deformations of $X$,
parametrized by an open set in $H^0(N_f)
\oplus H^1(W)$, and that 
this family yields deformations of $X$ as a Blanchard-Calabi 3-fold
(one can show that these are the deformations of $X$ which preserve the
fibration $\pi$).

All that remains is thus to show that the Kodaira Spencer map of this
family is surjective.  

Since we already remarked that we have an isomorphism
$W \cong \pi_* (\Theta_X) $, then $H^1(W) \cong H^1\pi_* (\Theta_X)$,
and we only have to control $H^0 (\mathcal{R}^1 \pi_*
\Theta_X )$.

To this purpose we split the long exact sequence of derived direct images
into the following exact sequences
$$0 \rightarrow (\Theta_C)\rightarrow V^{\vee} \otimes W \rightarrow
\mathcal{R} \rightarrow 0$$ 
$$0 \rightarrow \mathcal{R} \rightarrow \mathcal{R}^1 \pi_*
\Theta_X  \rightarrow
\mathcal{K} \rightarrow 0$$ 
$$0 \rightarrow \mathcal{K} \rightarrow ((\Theta_C)\otimes V^{\vee} ) 
\rightarrow
W \otimes \Lambda^2 (V^{\vee} )$$ 

and observe that all we need to prove is the vanishing 
$ H^0(\mathcal{K}) = 0.$

In fact, then $H^0(N_f) = H^0(\mathcal{R}) = H^0 (\mathcal{R}^1 \pi_*
\Theta_X )$ .

By the last exact sequence, it would suffice to show  the injectivity
of the linear map $H^0 ((\Theta_C)\otimes V^{\vee} ) \rightarrow
H^0(W \otimes \Lambda^2 (V^{\vee} )$.

But indeed we will show that the sheaf $\mathcal{K} $ is the zero sheaf.

Consider again in fact the beginning of the exact sequence of derived direct
images
$$ 0 \rightarrow (W) \rightarrow  \pi_* \Theta_X   \rightarrow 
(\Theta_C)  \rightarrow $$
$$  \rightarrow (V^{\vee} \otimes W):$$
and observe then that the homomorphism $(\Theta_C)  \rightarrow  (V^{\vee}
\otimes W)=
 \mathcal{R}^1 \pi_*
(\pi^*(W) )$ is indeed
the derivative of the map $f : C \rightarrow G$.

Since the
homomorphism of sheaves we are considering is 
$(\Theta_C)\otimes V^{\vee}  =   \mathcal{R}^1 \pi_*
(\pi^*(\Theta_C) ) \rightarrow  \mathcal{R}^2 \pi_*
(\pi^*(W) ) = 
W \otimes \Lambda^2 (V^{\vee} )$ the sheaf map $(\Theta_C)\otimes V^{\vee}
 \rightarrow W \otimes \Lambda^2 (V^{\vee} $is induced by wedge
product of $(\Theta_C)  \rightarrow  (V^{\vee}
\otimes W)$ with the identity of $V^{\vee}$, therefore we will actually
have that  the sheaf $\mathcal{K}$ is zero if the subsheaf of 
$ (V^{\vee} \otimes W)$  given by the image of $(\Theta_C)  \rightarrow 
(V^{\vee} \otimes W)$  consists of tensors of generical rank $2$.

But this holds by virtue of the hypothesis that $X$ be non developable.

Finally, let $X$ be a Sommese-Blanchard 3-fold with $L$ of degree $d \geq
g$ : then we have the desired vanishing  $H^1(V^{\vee} \otimes W) = 0$.

The condition that the associated ruled surface is non developable follows
from a direct calculation which shows that, if $f(t)$ is given as the
subspace generated by the  columns of a  $4 \times 2$ matrix $B(t)$,
then the $4 \times 4$  matrix $B(t) B'(t)$ has determinant equal to the
square (up to sign) of the Wronskian determinant of the section $s$.

But $s$ is non constant, for a Blanchard-Calabi 3-fold, therefore the
Wronskian determinant is not identically zero and we are done.

\hfill Q.E.D. \\

\begin{oss}
How does the developable case occur and which are its small deformations?
The condition that the four sections should in every point of $C$ provide a
real basis of the fibre of $W$ simply means that there is no real point in
the developable surface associated to $f: C \rightarrow G$.
\end{oss}

A more detailed analysis of the developable Blanchard-Calabi 3-folds and of
their deformations could allow  a positive answer to the following

QUESTION: Do the above Blanchard-Calabi 3-folds provide infinitely many
deformation types on the same differentiable manifold $C \times T$?
( In other words, do these just give infinitely many irreducible 
components of the "moduli space", or also "connected" components?)

 Calabi's construction is also quite beautiful, so that we cannot refrain
from indicating its main ideas.

The first crucial observation is that, interpreting $\R^7$ as the space
 $C^i$ of imaginary Cayley numbers, for any oriented hypersurface $M
\subset \R^7$, the  Cayley product produces an almost complex structure as
follows:

$$ J(v) = (v \times n)^i, $$ 

where $v$ is a tangent vector at the point  $x\in M$, $n$ is the normal 
vector at   $x\in M$, and $w^i$ stands for the imaginary part of a Cayley
number $w$.

For instance , this definition provides the well known non integrable
almost complex structure on the 6-sphere $S^6$.

Moreover, Calabi shows that the almost complex structure is special
hermitian, i.e., that its canonical bundle is trivial.

Second, Calabi analyses when is the given complex structure integrable,
proving in particular the following.  If we write $\R^7 = \R^3 \times
\R^4$, according to the decomposition $C = H \oplus He_5$ ( here $H$ is
the space of Hamilton's quaternions), and $M$ splits accordingly as a
product $S\times \R^4$, then the given complex structure is integrable if
and only if $S$ is a minimal surface in $\R^3$.

The third ingredient is now a classical method used by Schwarz in order
to construct minimal surfaces. Namely, let $C$ be a hyperelliptic curve,
so that the canonical map is a double cover of a rational normal curve
of degree $g-1$. If $g=3$ we have  then a basis of
holomorphic differentials 
$\omega_1, \omega_2,\omega_3$ such that the sum of their squares equals
zero. Similarly, for every odd genus, via an appropriate projection  to
$\PP ^2$, we obtain three linearly independent holomorphic differentials,
without common zeros, and satisfying also the relation $\omega_1^2 +
\omega_2^2 +\omega_3^2 = 0$.

The integral of the real parts of the $\omega_i$'s provides a multivalued
map of  $C$ to $\R^3$ which is a local embedding. Since another local
determination differs by translation, the tangent space to a point in
$C$ is then naturally a subspace of $\R^3$, whence Cayley multiplication
provides a well defined complex structure on $C \times \R^4$.
Since moreover this complex structure is also invariant by translations
on $\R^4$, we can descend a complex structure on $C\times T$.

From  the contructions of Blanchard and Calabi we deduce a
negative answer to the problem mentioned in remark 2.1 and in the
introduction.

\begin{cor}
Products of a curve of genus $g \geq 1$ with a 2-dimensional complex
torus provide  examples of  complex manifolds which are a $K(\pi, 1)$,
and for which there are different deformation types.
\end{cor}

We shall see in the next section that the answer continues to be negative
even if we restrict to K\"ahler, and indeed projective manifolds, even in
dimension = 2.

\section{Moduli spaces of surfaces of general type}
\label{sixth}
In this section we  begin to describe a recent result ([Cat8]), showing
the existence of  complex surfaces which are $K(\pi, 1)$'s and for which
there are different deformation types.

Before we get into the details of the construction, it seems
appropriate to give a more general view of the status of the art
concerning deformation, differentiable and topological types of
algebraic surfaces of general type.

Let $S$ be a minimal surface of general type: then  to $S$  we attach two
positive integers $ x = \chi (\hol_S)$,
$ y =K^2_S$ which are invariants of the oriented topological type
of $S$. 

The moduli space of the surfaces with invariants $(x,y)$ is a
quasi-projective variety (cf. \cite{gie}) defined over the integers, in
particular it is a real variety.

Fixed $(x,y)$ we have several possible topological types,
but indeed only two if moreover the surface $S$ is simply
connected. These two cases are distinguished as follows:

\begin{itemize}
\item
 $S$ is EVEN, i.e., its intersection form is even :
then $S$ is a connected sum of copies of a $K3$ surface
and of $ \PP ^1_{\C} \times \PP ^1_{\C}$.
\item
$S$ is ODD : then $S$ is a connected sum of copies of
$\PP ^2_{\C} $ and ${\PP ^2_{\C}}^{opp} .$
\end{itemize}

\begin{oss}
${\PP ^2_{\C}}^{opp} $ stands for the same manifold as
${\PP ^2_{\C}}$, but with reversed orientation.

It is rather confusing, especially if one has to do with real structures,
that some authors use the symbol
$\bar{\PP ^2_{\C}} $ for  ${\PP ^2_{\C}}^{opp} $ .
\end{oss}

In general, the fundamental group is a powerful topological
invariant. Invariants of the differentiable structure have been
found by Donaldson and Seiberg -Witten, and one can easily show
that on a connected component of the moduli space the
differentiable structure remains fixed.

Up to recently, the converse question DEF = DIFF ? was open,
but recently counterexamples have been given, by Manetti (\cite{man3}) for
simply connected surfaces, by Kharlamov and Kulikov (\cite{k-k}) for rigid
surfaces, while we have found rather simple examples (\cite{cat8}):

\begin{teo}
Let $S$ be a surface isogenous to a product, i.e., a quotient 
$ S = (C_1 \times C_2) /G$ of a product of curves  by the
free action of a finite group $G$. Then any surface with the
same fundamental group as $S$ and the same Euler number of $S$
is diffeomorphic to $S$. The corresponding moduli space 
$ M^{top}_S$ is either irreducible and connected or it contains
two connected components which are exchanged by complex
conjugation. There are infinitely many examples of the latter
case.
\end{teo}

\begin{cor}
1) DEF $\neq$ DIFF.

2) There are moduli spaces without real points 

3) There are complex surfaces whose  fundamental
group cannot be the fundamental group of a real surface.
\end{cor}

For the construction of these examples 
we imitate the hyperelliptic surfaces, in the sense that we
take $ S = (C_1 \times C_2) /G$ where $G$ acts freely on $C_1$,
whereas the quotient $ C_2 /G$ is $\PP^1_{\C}$.

Moreover, we assume that the projection $\phi: C_2 \rightarrow
\PP^1_{\C}$ is branched in only three points, namely, we have a
so called TRIANGLE CURVE.

What happens is that if two surfaces of such sort are
antiholomorphic, then there would be an antiholomorphism of the
second triangle curve (which is rigid).

Now, giving such a branched cover $\phi$ amounts to viewing
the group $G$ as a quotient of the free group on two elements.
Let $a, c$ be the images of the two generators, and set 
$ abc = 1$. 

We  find such a $G$ with the properties that the respective
orders of $a,b,c$ are distinct, whence an antiholomorphism of
the triangle curve would be a lift of the standard complex
conjugation if the 3 branch points are chosen to be real, e.g.
$-1, 0 -+1$. 

Such a lifting exists if and only if the group $G$ admits
an automorphism $\tau$ such that $\tau (a) = a^{-1}, \tau (c) =
c^{-1}$. 

An appropriate semidirect product will be the desired group for which
such a lifting does not exist.

For this reason, in the next section we briefly recall the notion and the
simplest examples of the so called triangle curves.

\section{Some non real triangle curves}
\label{sixth}

Consider  the set $B \subset \PP^1_{\C}$ consisting of  three 
real points $ B : = \{-1, 0, 1\}.$ 

We choose $2$ as a base point in $\PP^1_{\C} - B$, and we take 
the following generators $ \alpha, \beta, \gamma$ of 
$ \pi_1 (\PP^1_{\C} - B, 2)$ such that

$$ \alpha \beta \gamma = 1 $$ 

and $ \alpha,  \gamma$ are  free generators 
of $ \pi_1 (\PP^1_{\C} - B, 2)$ with  $ \alpha,  \beta$ as indicated in the
following picture

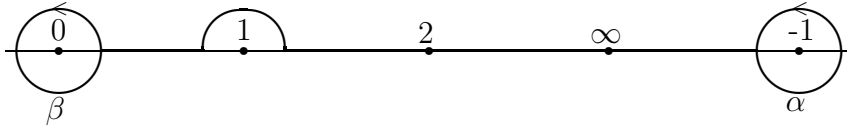
\font\mate=cmmi8
\def\min{\hbox{\mate \char60}\hskip1pt}
\begin{center}
\begin{picture}(350,100)(0,-40)
%
%
\put(0,0){\line(1,0){300}}
\put(238,0){\line(1,0){80}}
%
%
\put(20,0){\circle*{3}}
\put(17,4){0}
\put(90,0){\circle*{3}}
\put(87,4){1}
\put(160,0){\circle*{3}}
\put(156,3){2}
\put(228,0){\circle*{3}}
\put(221,3){$\infty$}
\put(300,0){\circle*{3}}
\put(296,4){-1}
{\thicklines
%
%
\put(20,0){\circle{30}}
\put(36,0){\line(1,0){39}}
\put(90,0){\oval(31,31)[t]}
\put(105,0){\line(1,0){55}}
\put(17,14){$\min$}
\put(15,-26){$\beta$}
%
%
\put(300,0){\circle{30}}
\put(160,0){\line(1,0){124}}
\put(297,14){$\min$}
\put(295,-23){$\alpha$}
}
\end{picture}
\end{center}

With this choice of basis , we have provided an isomorphism
of $ \pi_1 (\PP^1_{\C} - B, 2)$ with the group
$$  T_{\infty} : = < \alpha, \beta, \gamma | \ 
 \alpha \beta \gamma = 1 >.$$
For each finite group $G$ generated by two elements $a,b$, 
passing from Greek to latin letters we obtain a tautological
surjection 

$ \pi : T_{\infty}  \rightarrow G .$ 

I.e., we set $\pi (\alpha) = a, \pi (\beta) = b$ and we define 
$ \pi (\gamma) := c.$ (then $abc=1$).
To $\pi$ 
we associate the Galois covering $f:  C \rightarrow 
\PP^1_{\C}$,  branched on $B$  and with group $G$.

Notice that the Fermat curve $ C := \{(x_0,x_1,x_2) \in \PP^2_{\C} |
x_0^n + x_1^n + x_2^n = 0\}$ is in two ways a triangle curve,
since we can take the quotient of $C$ by the group $G := (\Z/n) ^2$ of 
diagonal projectivities with entries n-th roots of unity, but
also by the full group $A= Aut(C)$ of automorphisms, which
is a semidirect product of the normal subgroup $G$ by the symmetric
group permuting the three coordinates.
For $G$ the three branching multiplicities are all equal to $n$,
whereas for $A$ they are equal to $(2,3,2n)$.

Another interesting example is provided by the {\bf Accola curve}
( cf. \cite{acc1}, \cite{acc2}),
the curve $ Y_g $ birational to the affine curve of equation

$$ y^2 = x^{2g+2} - 1 .$$

If we take the group $G \cong \Z/ 2 \times \Z/(2g+2)$ which acts 
multiplying $y$ by $-1$, respectively $x$ by a primitive $2g+2$-root
of $1$, we realize $ Y_g $ as a triangle curve with
 branching multiplicities $ (2, 2g+2, 2g+2)$. $G$ is not however
the full automorphism group, in fact if we add the involution
sending $x$ to $1/x$ and $y$ to $iy/x^{g+1}$, then we get the direct
product $\Z/ 2 \times D_{2g+2}$ (which is indeed the full group
of automorphisms of $ Y_g $ as it is well known and as also 
follows from the
next lemma), a group which represents $ Y_g $ as a triangle curve with
 branching multiplicities $( 2, 4 , 2g+2)$.

One can get many more examples by taking (a la Macbeath. \cite{mcb})
unramified coverings of the above curves (associated to characteristic
subgroups of the fundamental group).
The following natural question arises then:
which are the curves which admit more than one realization as
triangle curves?
It is funny to observe:
\begin{lem}
Let $f : C \rightarrow \PP^1_{\C} = C/G$ be a triangle covering
where the branching multiplicities $m,n,p$ are all distinct 
(thus we assume $m < n < p$).
Then the group $G$ equals the full group $A$ 
of automorphisms of $C$.
\end{lem}
 Idea:

I. By Hurwitz's formula the cardinality of $G$ is in general
given by the formula $$ |G| = 2 (g-1) (1-1/m-1/n-1/p)^{-1} .$$

II. Assume that $A\neq G$ and let $F: \PP^1_{\C} = C/G
\rightarrow \PP^1_{\C} = C/A$ be the induced map.
Then $f': C \rightarrow \PP^1_{\C} = C/A$ is again a triangle covering,
otherwise the number of branch points would be $\geq 4$ and
we would have a non trivial family of such Galois covers
with group $A$.

III. We claim now that the three branch points of $f$ 
cannot have distinct images
through $F$: otherwise the branching multiplicities  $m' \leq n'
\leq p'$ for $f'$ would be not less than the respective multiplicities
for $f$, and by the analogous of formula I for $|A|$ we would
obtain $ |A| \leq |G|$, a contradiction.
The rest of the proof is complicated.
\hfill End.

We come now to our particular triangle curves.
Let $r, \ m$ be positive integers  $r\geq 3 , m \geq 4$ and set 
$$ p := r^m -1,\  n:= (r-1) m \ . $$
Notice that the three integers $m < n < p $ are distinct.
Let $G$ be the following semidirect product of $ \Z/p $ by $ \Z/m$:

$$  G : = < a,c \ |\  a^m = 1 , \ c^p = 1, \ a c a^{-1} = c^r \ > $$ 
One sees easily that the period of
$b$ is exactly $n$. 

\begin{prop}
The triangle curve $ C$ associated to $\pi$ is not antiholomorphically
 equivalent to itself ( i.e., it is not isomorphic to its conjugate). 
\end{prop}
Idea:

We derive a contradiction assuming the existence of
  an antiholomorphic automorphism $\sigma$ of $C$.
  
STEP I 
: $G$  = $A$, where $A$ is the group 
of holomorphic automorphisms of $C$, $A: = Bihol(C,C)$.

\Proof This follows from the previous lemma 2.3. 

STEP II : if $\sigma$ exists, it must be a lift of complex
conjugation.

\Proof In fact $\sigma$ normalizes $Aut(C)$, whence it must induce
a antiholomorphism of $\PP^1_{\C}$ which is the identity on
$B$, and therefore must be complex conjugation.

STEP III : complex conjugation does not lift.

\Proof This is purely an argument about covering spaces:
complex conjugation acts on $ \pi_1 (\PP^1_{\C} - B, 2) \cong 
T_{\infty} $,
as it is immediate to see with our choice of basis,
by  the automorphism $\tau$ sending $\alpha$,  $\gamma$ 
to their respective inverses.

Thus, complex conjugation lifts if and only if $\tau$ preserves
the normal subgroup $K : = ker \pi$. In turn, this means that there
is an automorphism $ \rho : G \rightarrow G$ with

$$ \rho (a) = a^{-1}, \  \rho (c) = c^{-1}. $$ 

Recall now the relation $ a c a ^{-1} = c^r $ : 
applying $\rho $,
we would get $ a^{-1} c ^{-1}a  = c^{-r }$, or, equivalently,

$$ a^{-1} c a  = c^{r }.$$

But then we would get $ a c a ^{-1} = c^r  =a^{-1} c a  =a^{m-1}
 c a^{m-1} = c^{r^{m-1} } $ , which holds only if

$$ r \equiv r^{m-1} \ (mod p).$$

 Since $p = r^m-1 $ we obtain,
 after multiplication by $r$,  that we should have $r^2 
\equiv  1 \ (mod \ p)$ but this is the desired  contradiction,
because $ r^2 - 1 < r^m - 1 = p $.

\qed

\section{The examples of surfaces isogenous to a product}
\label{seventh}
\begin{df}
A projective surface $S$ is said to be isogenous to a  (higher) product 
if it admits a finite unramified covering by
 a product of curves of genus $\geq 2$. 
 \end{df}
 
 \begin{oss}
 In this case, (cf. \cite{cat7}, props.
 3.11 and 3.13) there
exist  Galois realizations $ S = (C_1 \times C_2 )/ G$, and  each
such Galois realization dominates a uniquely determined minimal
one. $S$ is said to be of nonmixed type if $G$ acts via a product
action of two respective actions on $C_1, C_2$. Else, $G$ contains a
subgroup $G^0$ of index $2$ such that $(C_1 \times C_2 )/ G^0$ is of
nonmixed type.
\end{oss}

\begin{prop}
Let $S, S'$ be surfaces isogenous to a higher product, and let
$\sigma : S \rightarrow S'$ be a antiholomorphic isomorphism.
Let moreover $ S = (C_1 \times C_2 )/ G$, 
$ S' = (C'_1 \times C'_2 )/ G'$ be the respective
minimal Galois realizations. Then, up to possibly exchanging
$C'_1 $ with $ C'_2 $, there exist antiholomorphic isomorphisms
$\tilde{\sigma}_i , \ i= 1,2 $ such that $\tilde{\sigma} :=
\tilde{\sigma}_1 \times \tilde{\sigma}_2  $ normalizes the action
of $G$, in particular $\tilde{\sigma}_i $ normalizes the action
of $ G^0 $ on $C_i$.
\end{prop}
\proof
Let us view $\sigma$ as yielding a complex isomorphism 
$\sigma : S \rightarrow \bar{S'}$. 
Consider the exact sequence corresponding  to the minimal
Galois realization 
$ S = (C_1 \times C_2 )/ G$,
$$1 \rightarrow H : = \Pi_{g_1} \times \Pi_{g_2} 
\rightarrow  \pi_1(S) \rightarrow G \rightarrow 1. $$ 
Applying $\sigma_*$  to it, we infer by theorem 3.4 of (\cite{cat7})
that we obtain an exact sequence associated to a Galois 
realization of $ \bar{S'}$. Since $\sigma$ is an isomorphism,
we get a minimal one, which is however unique.
Whence,we get an isomorphism    $ \tilde{\sigma} :
(C_1 \times C_2 )  \rightarrow (\bar{C'_1} \times \bar{C'_2} ) $,
which is of product type by the rigidity lemma (e.g., lemma 3.8
of \cite{cat7}). Moreover this isomorphism must normalize
the action of $ G \cong G'$, which is exactly what we claim.

\qed

In (\cite{cat7}, cf. correction in \cite{cat8}) we have proven:

\begin{teo}
Let $S$ be a surface isogenous to a product, i.e., a quotient 
$ S = (C_1 \times C_2) /G$ of a product of curves  by the
free action of a finite group $G$. Then any surface $S'$ with the
same fundamental group as $S$ and the same Euler number of $S$
is diffeomorphic to $S$. The corresponding moduli space 
$ M^{top}_S = M^{diff}_S $ is either irreducible and connected 
or it contains
two connected components which are exchanged by complex
conjugation. 
\end{teo}

We are now going to explain  the construction
of our examples :

Let $G$ be the semidirect product group we constructed in section
$2$, and let $C_2$ be the corresponding triangle curve.

Let moreover $g'_1$ be any number greater or equal to $2$,
and consider the canonical epimorphism $\psi$ of $\Pi_{g'_1}$ 
onto a free group of rank $g'_1$, such that in terms of the standard 
bases $ \ a_1 , b_1, .... a_{g'_1},\ b_{g'_1} $, respectively  
$ \gamma_1,.. \gamma_{g'_1}$,
we have 
$\psi (a_i) = \psi(b_i) = \gamma_i$.

Compose then $\psi$  with any epimorphism of the free group
onto $G$, e.g. it suffices to compose with any $\mu$ such that 
$ \mu ( \gamma_1) = a$,  $\mu ( \gamma_2) = b$ ( and 
$\mu ( \gamma_j)$ can be chosen whatever we want for $j \geq 3$).

For any point $C'_1$ in the Teichm\"uller space we obtain a canonical
covering associated to the kernel of the epimorphism $ \mu \circ 
\psi : \Pi_{g'_1} \rightarrow G$, call it $C_1$.

\begin{df}
Let $S$ be the surface $S : =(C_1 \times C_2) / G$ ($S$ is smooth
because $G$ acts freely on the first factor).
\end{df}

\begin{teo}
For any two choices  $C'_1(I) , C'_1(II) $ 
of $C'_1$ in the Teichm\"uller space we get surfaces  $S(I), S(II)$
such that $S(I)$ is never isomorphic to $\bar{S}(II)$.
Varying $C'_1$ we get a connected component of the moduli space,
which has only one other connected component, given by the conjugate of
the previous one.
\end{teo}

The last result that we have obtained as an application of these ideas is
the following puzzling:

\begin{teo}
Let $S$ be a surface in one of the families constructed above.
Assume moreover that $X$ is another complex surface such
that $\pi_1(X) \cong \pi_1(S)$. Then $X$ does not admit any real
structure.
\end{teo}
\proof(idea)

Observe that since $S$ is a classifying space for the fundamental
group of $\pi_1(S)$, then by the isotropic subspace theorem of
(\cite{cat4})
 the Albanese mapping of $X$ maps onto a curve $C'(I)_2$
of the same genus as $C'_2$.

Consider now the unramified covering $\tilde{X}$ associated to
the kernel of the epimorphism $\pi_1(X) \cong \pi_1(S) \rightarrow G$.

Again by the isotropic subspace theorem, there exists a holomorphic
map $\tilde{X} \rightarrow C(I)_1 \times C(I)_2$, where moreover
the action of $G$ on $\tilde{X}$ induces  actions of $G$ on
both factors which either have the same topological types as the actions
of $G$ on $C_1$, resp. $C_2$ , or have both the topological types
of the actions on the respective conjugate curves.

By the rigidity of the triangle curve $C_2$, in the former case
$ C(I)_2 \cong  C_2$, in the latter $ C(I)_2 \cong  \bar{C_2}$.

Assume  now that $X$ has a real structure $\sigma$: then the same
argument as 
in \cite{c-f2}  section 2 shows that $\sigma$ induces  a product
antiholomorphic map $\tilde{\sigma} : C(I)_1 \times C(I)_2 
\rightarrow C(I)_1 \times C(I)_2$. In particular, we get a 
non costant antiholomorphic map of $C_2$ to itself, 
contradicting proposition 2.3 .
 \qed 

In \cite{cat7}, \cite{cat8} we gave the following

\begin{df}
A Beauville surface is a rigid surface isogenous to a product.
\end{df}

Beauville gave examples of these surfaces (\cite{bea}), as quotients of the
product $C\times C$ where $C$ is the Fermat quintic curve. This example
is real.

It would be interesting to

\begin{itemize}
\item
Classify all the Beauville surfaces, at least those of non mixed type.
\item
Classify all the non real Beauville surfaces, especially those which are
not isomorphic to their conjugate surface.

\end{itemize}

\bigskip

{\bf Acknowledgements.}
 
I am grateful to  S.T. Yau for pointing out Calabi's example, and for
his invitation to Harvard, where some of the results in this paper were
obtained.

Thanks also to P. Frediani, with whom I  started to investigate the
subtleties of the "real" world ([C-F1] and [C-F2]), for several
interesting conversations.

I would like  moreover to thank Mihai Paun for pointing out a small error, 
and especially express my  gratitude to Andrew Sommese for reminding me
about the relevance of Blanchard's examples.

\vfill

\noindent
{\bf Author's address:}

\bigskip

\noindent 
Prof. Fabrizio Catanese\\
Lehrstuhl Mathematik VIII\\
Universit\"at Bayreuth\\
 D-95440, BAYREUTH, Germany

e-mail: Fabrizio.Catanese@uni-bayreuth.de

\end{document}